
\def\LaTeX{%
  \let\Begin\begin
  \let\End\end
  \let\salta\relax
  \let\finqui\relax
  \let\futuro\relax}

\def\UK{\def\our{our}\let\sz s}
\def\USA{\def\our{or}\let\sz z}



\LaTeX

\USA


\salta

\documentclass[twoside,12pt]{article}
\setlength{\textheight}{24cm}
\setlength{\textwidth}{16cm}
\setlength{\oddsidemargin}{2mm}
\setlength{\evensidemargin}{2mm}
\setlength{\topmargin}{-15mm}
\parskip2mm


\usepackage{amsmath}
\usepackage{amsthm}
\usepackage{amssymb}
\usepackage[mathcal]{euscript}

\usepackage[usenames,dvipsnames]{color}
%
%
%

\def\pier{\color{red}}
\def\js{\color{blue}}
%
%

\let\pier\relax
\let\js\relax


\bibliographystyle{plain}


%

\finqui

\def\Beq{\Begin{equation}}
\def\Eeq{\End{equation}}
\def\Bsist{\Begin{eqnarray}}
\def\Esist{\End{eqnarray}}

\def\Bthm{\Begin{theorem}}
\def\Ethm{\End{theorem}}
\def\Blem{\Begin{lemma}}
\def\Elem{\End{lemma}}
\def\Bprop{\Begin{proposition}}
\def\Eprop{\End{proposition}}

\def\Brem{\Begin{remark}\rm}
\def\Erem{\End{remark}}

\def\Bdim{\Begin{proof}}
\def\Edim{\End{proof}}
\let\non\nonumber




\def\step #1 \par{\medskip\noindent{\bf #1.}\quad}
\def\proofstep #1 \par{\medskip\noindent{\emph{#1.}}\ }


\def\Lip{Lip\-schitz}
\def\holder{H\"older}
\def\aand{\quad\hbox{and}\quad}

\def\lsc{lower semicontinuous}

\def\lhs{left-hand side}
\def\rhs{right-hand side}
\def\sfw{straightforward}
\def\wk{well-known}


\def\organiz{organi\sz}

\def\regulariz{regulari\sz}


\def\multibold #1{\def\arg{#1}%
  \ifx\arg\pto \let\next\relax
  \else
  \def\next{\expandafter
    \def\csname #1#1#1\endcsname{{\bf #1}}%
    \multibold}%
  \fi \next}

\def\pto{.}

\def\multical #1{\def\arg{#1}%
  \ifx\arg\pto \let\next\relax
  \else
  \def\next{\expandafter
    \def\csname cal#1\endcsname{{\cal #1}}%
    \multical}%
  \fi \next}


\def\multimathop #1 {\def\arg{#1}%
  \ifx\arg\pto \let\next\relax
  \else
  \def\next{\expandafter
    \def\csname #1\endcsname{\mathop{\rm #1}\nolimits}%
    \multimathop}%
  \fi \next}

\multibold
qwertyuiopasdfghjklzxcvbnmQWERTYUIOPASDFGHJKLZXCVBNM.

\multical
QWERTYUIOPASDFGHJKLZXCVBNM.

\multimathop
dist div dom meas sign supp .


\def\accorpa #1#2{\eqref{#1}--\eqref{#2}}
\def\Accorpa #1#2 #3 {\gdef #1{\eqref{#2}--\eqref{#3}}%
  \wlog{}\wlog{\string #1 -> #2 - #3}\wlog{}}


\def\QED{\hfill $\square$}

\def\tonde #1{\left(#1\right)}

\def\graffe #1{\mathopen\{#1\mathclose\}}

\def\<#1>{\mathopen\langle #1\mathclose\rangle}
\def\norma #1{\mathopen \| #1\mathclose \|}

\def\iot {\int_0^t}
\def\ioT {\int_0^T}
\def\iO{\int_\Omega}
\def\intQt{\int_{Q_t}}
\def\intQ{\int_Q}

\def\dt{\partial_t}
\def\dn{\partial_\nu}
\def\ds{\,ds}

\def\cpto{\,\cdot\,}

\def\checkmmode #1{\relax\ifmmode\hbox{#1}\else{#1}\fi}
\def\aeO{\checkmmode{a.e.\ in~$\Omega$}}
\def\aeQ{\checkmmode{a.e.\ in~$Q$}}

\def\aat{\checkmmode{for a.a.~$t\in(0,T)$}}


\def\erre{{\mathbb{R}}}




\def\genspazio #1#2#3#4#5{#1^{#2}(#5,#4;#3)}
\def\spazio #1#2#3{\genspazio {#1}{#2}{#3}T0}

\def\L {\spazio L}
\def\H {\spazio H}
\def\W {\spazio W}

\def\C #1#2{C^{#1}([0,T];#2)}

\def\Vp{V^*}


\def\Lx #1{L^{#1}(\Omega)}
\def\Hx #1{H^{#1}(\Omega)}
\def\Wx #1{W^{#1}(\Omega)}
\def\Cx #1{C^{#1}(\overline\Omega)}

\def\Luno{\Lx 1}
\def\Ldue{\Lx 2}
\def\Linfty{\Lx\infty}
\def\Lq{\Lx4}
\def\Huno{\Hx 1}
\def\Hdue{\Hx 2}


\def\LQ #1{L^{#1}(Q)}
\def\CQ {C^0(\overline Q)}


\let\theta\vartheta
\let\eps\varepsilon

\let\TeXchi\chi                         
\newbox\chibox
\setbox0 \hbox{\mathsurround0pt $\TeXchi$}
\setbox\chibox \hbox{\raise\dp0 \box 0 }
\def\chi{\copy\chibox}


\def\muz{\mu_0}
\def\rhoz{\rho_0}

\def\normaV #1{\norma{#1}_V}
\def\normaH #1{\norma{#1}_H}
\def\normaW #1{\norma{#1}_W}

\def\rmin{r_*}
\def\smin{s_*}
\def\Kz{K_*}

\def\kamu{\kappa(\mu)}
\def\kmin{\kappa_*}
\def\kmax{\kappa^*}
\def\rmin{r_*}

\def\coeff{1+2g(\rho)}
\def\coeffr{1+2g(r)}
\def\coeffs{1+2g(\rho(s))}
\def\coefft{1+2g(\rho(t))}

\def\T{\calT_\tau}

\def\mut{\mu_\tau}
\def\rhot{\rho_\tau}
\def\xit{\xi_\tau}
\def\kappat{\kappa_\tau}
\def\Kt{K_\tau}
\def\kamut{\kappat(\mut)}

\def\coefftau{1+2g(\rhot)}

\def\coeffn{1+2g(\rhon)}
\def\coeffnt{1+2g(\rhon(t))}

\def\mun{\mu_n}
\def\rhon{\rho_n}
\def\xin{\xi_n}
\def\munmu{\mu_{n-1}}

\def\coeffn{1+2g(\rhon)}

\def\iotmt{\int_0^{t-\tau}}
\def\intQtmt{\iotmt\!\!\!\iO}



\DeclareMathAlphabet{\mathbf}{OT1}{cmr}{bx}{it}

\newcommand{\hb}{\mathbf{h}}
\newcommand{\nb}{\mathbf{n}}
\newcommand{\Mb}{\mathbf{M}}
\newcommand{\Hb}{\mathbf{H}}
 \font\mba=cmmib10 scaled
\magstephalf

\def\csib{\hbox{\mba {\char 24}}}
\Begin{document}



\title{\bf 
Global existence and uniqueness\\
for a singular/degenerate\\
Cahn-Hilliard system with viscosity%
}
\author{}
\date{}
\maketitle
\begin{center}
\vskip-2cm
{\large\bf Pierluigi Colli$^{(1)}$}\\
{\normalsize e-mail: {\tt pierluigi.colli@unipv.it}}\\[.25cm]
{\large\bf Gianni Gilardi$^{(1)}$}\\
{\normalsize e-mail: {\tt gianni.gilardi@unipv.it}}\\[.25cm]
{\large\bf Paolo Podio-Guidugli$^{(2)}$}\\
{\normalsize e-mail: {\tt ppg@uniroma2.it}}\\[.25cm]
{\large\bf J\"urgen Sprekels$^{(3)}$}\\
{\normalsize e-mail: {\tt sprekels@wias-berlin.de}}\\[.45cm]
$^{(1)}$
{\small Dipartimento di Matematica ``F. Casorati'', Universit\`a di Pavia}\\
{\small via Ferrata 1, 27100 Pavia, Italy}\\[.2cm]
$^{(2)}$
{\small Dipartimento di Ingegneria Civile, Universit\`a di Roma ``Tor Vergata''}\\
{\small via del Politecnico 1, 00133 Roma, Italy}\\[.2cm]
$^{(3)}$
{\small WIAS -- Weierstra\ss-Institut f\"ur Angewandte Analysis und Stochastik}\\
{\small Mohrenstra\ss e\ 39, 10117 Berlin, Germany}\\[.8cm]
\end{center}


\Begin{abstract}
Existence and uniqueness are investigated for a nonlinear diffusion problem of 
phase-field type, consisting of a parabolic system of two partial differential 
equations, complemented by Neumann homogeneous boundary conditions and initial 
conditions. 
This system aims to model two-species phase segregation on an atomic 
lattice~\cite{Podio}; in the balance equations of microforces and microenergy, the two 
unknowns are the order parameter $\rho$ and the chemical potential $\mu$. A 
simpler version of the same system has {\js recently been} discussed in \cite{CGPS3}.
In this paper, a fairly more general phase-field equation for $\rho$ is coupled with a genuinely 
nonlinear diffusion equation for $\mu$. The existence of a global-in-time solution 
is proved with the help of suitable a priori estimates. In the case of 
a constant atom mobility, a new and rather unusual uniqueness 
proof is given, based on a suitable combination of variables.
\End{abstract}


\salta

\pagestyle{myheadings}
\newcommand\testopari{\sc Colli \ --- \ Gilardi \ --- \ Podio-Guidugli \ --- \ Sprekels}
\newcommand\testodispari{\sc Existence -- uniqueness for a Cahn-Hilliard system  
with viscosity}
\markboth{\testodispari}{\testopari}

\finqui


\section{Introduction}
\label{Intro}
\setcounter{equation}{0}

In this paper, the last so far in a series \cite{CGPS1,CGPS2,CGPS3,CGPS4,CGPS5,CGPS7},  we further our mathematical analysis of a mechanical model proposed by one of us \cite{Podio} for phase segregation through atom rearrangement on a lattice. On postponing a detailed presentation of the model and its antecedents until next section, we begin by pointing out what features of the system we study are more general, and therefore more difficult to handle mathematically, than in our previous paper \cite{CGPS3}. 

%
%
The initial and boundary value problem we here tackle consists in looking 
for two \emph{fields},  the {\it chemical potential} $\mu>0$ and the \emph{order 
parameter} $\rho\in (0,1)$, solving 
\Bsist
  \hskip-1cm&2 h(\rho) \,  \dt\mu
  + \mu \, h'(\rho) \, \dt\rho
  - \div \bigl( \kamu\nabla\mu \bigr) = 0
  & \quad \hbox{in $\Omega \times (0,T)$,}
  \label{Iprima}
  \\
   \hskip-1cm&\delta \dt\rho - \Delta\rho + f'(\rho)= \mu h'(\rho), 
   & \quad \hbox{in $\Omega \times (0,T)$,}
   \label{Iseconda}
  \\
  \hskip-1cm& (\kamu\nabla\mu)\cdot\nu|_\Gamma = \dn\rho|_\Gamma = 0
  & \quad \hbox{on $\Gamma \times (0,T)$,}
  \label{Ibc}
  \\
  \hskip-1cm& \mu(\cpto,0) = \muz
  \aand
  \rho(\cpto,0) = \rhoz
  & \quad \hbox{in $\Omega$,}
  \label{Icauchy}
\Esist
\Accorpa\Ipblnew Iprima Icauchy
where $\Omega$ denotes a bounded domain of $ \erre^3$ with conveniently smooth boundary $\Gamma$, $T>0$, and $\dn$ denotes differentiation in the direction of the outward normal $\nu$. 
In \eqref{Iprima}, the \emph{atom mobility} is specified by a nonnegative, continuous 
and  bounded, nonlinear function $\kappa$ of $\mu $ (in particular, the degeneracy of $\kappa$ around the critical value $\mu=0$ is admitted). The problem is parameterized by two scalar-valued functions, $h$ and $f$, and two positive numbers, $\varepsilon$ and $\delta$, both intended to be small. The parameter functions enter into the definition of the system's \emph{free energy}
\begin{equation}
\label{fe-2}
\psi=\widehat\psi(\rho,\nabla\rho,\mu) 
= - \mu \, h(\rho) + f (\rho) +\frac{1}{2}|\nabla\rho|^2,
\end{equation}
where the last two terms favor phase segregation, the former because it introduces local energy minima and the latter because it penalizes spatial changes of the order parameter (we have set equal to 1 the relative material constant). For
 $h$, one can take any smooth function provided it is bounded from below by a positive constant:
\begin{equation}\label{newh}
h(\rho) \geq \displaystyle\frac{\eps}2; 
\end{equation}
%
%
for $f$, the sum
\[
f(\rho)=f_1(\rho) + f_2(\rho) 
\] 
of a convex and \lsc\ function $f_1$, with proper domain $D(f_1) {\pier{}\subseteq {}}\erre $, and of a smooth function $f_2$ with no convexity properties, so as to allow for a double or multi-well potential $f$. Note that $f_1$ need not be differentiable in its domain, so that its possibly multivalued subdifferential  $\beta:=\partial f_1$ may appear in \eqref{Iseconda} in place of $f'_1$; in general, $\beta$ is only a graph, 
not necessarily a function, and it may include vertical (and horizontal) lines, 
as for example when
\begin{equation}
\label{ex2}
f_1 (\rho) =  I_{[0,1]} (\rho) =  \ \left\{
\begin{array}{ll}
0 & \text{if \ $0\leq \rho \leq 1$}
\\[0.1cm]
+\infty \ &\text{elsewhere}
\end{array}
\right. 
\end{equation}
{\pier and $\beta = \partial I_{[0, 1]}$ is specified by
\begin{equation}
\label{ex0}
\xi \in \beta (\rho) 
\quad \hbox{ if and only if } \quad 
\xi \ \left\{
\begin{array}{ll}
\displaystyle
\leq \, 0 \   &\hbox{if } \ \rho=0   
\\[0.1cm]
= \, 0 \   &\hbox{if } \ 0< \rho < 1  
\\[0.1cm]
\geq \, 0 \  &\hbox{if } \  \rho  = 1  
\\[0.1cm]
\end{array}
\right. .
\end{equation}
} 

The simpler situation dealt with in \cite{CGPS3} obtains for $\kappa$ constant-valued (and hence set equal to 1, without any loss of generality),
\begin{equation}\label{oldh}
h(\rho)=\rho,
\end{equation}
and $f$ a double-well potential defined in $(0,1)$, whose derivative $f'$ is singular at the endpoints $\rho=0$ and $\rho=1$: e.g.,
\begin{equation}
\label{ex1} 
f(\rho)= \alpha_1 \, \{ \rho\,\ln (\rho)+(1-\rho)\,\ln (1-\rho) \} + \alpha_2 \, \rho \, (1-\rho)
\end{equation}
for some positive constants $\alpha_1$ and $\alpha_2$.\footnote{Note that, according to whether or not $\alpha_1 \geq 2\alpha_2$, it turns out that $f$ is convex in the whole of $[0,1]$ or it exhibits two wells with a local maximum at $\rho=1/2$.}
Under these less general circumstances, system \Ipblnew\ reduces to
\Bsist
  & \eps\, \dt\mu + 2\rho \,\dt\mu + \mu \, \dt\rho - \Delta\mu = 0
  & \quad \hbox{in $\Omega \times (0,T)$,}
  \label{Iprima-3}
  \\
  & \delta\, \dt\rho - \Delta\rho + f'(\rho) = \mu 
  & \quad \hbox{in $\Omega \times (0,T)$,}
  \label{Iseconda-3}
  \\
  & \dn\mu = \dn\rho = 0
  & \quad \hbox{on $\Gamma \times (0,T)$,}
  \label{Ibc-3}
  \\
  & \mu(\cpto,0) = \muz
  \aand
  \rho(\cpto,0) = \rhoz
  & \quad \hbox{in $\Omega$.}
  \label{Icauchy-3}
\Esist
\Accorpa\Ipbl Iprima-3 Icauchy-3
Note that $h$ might attain its lower 
bound  for some significant values of $\rho$, that is, 
for some $\rho$'s lying in the domain of $f_1$: actually, this was the case for 
$h$ defined as in \eqref{oldh} over the interval $[0,1]$, that is, over 
the effective domain of both potentials in \eqref{ex1} and 
\eqref{ex2}. We were prompted to generalize \eqref{oldh} as in \eqref{newh} by an interesting remark of Alexander Mielke, when one of us was 
lecturing on our results, namely, that the behavior of 
$$ h(\rho) = \rho + \hbox{ small parameter} $$
is different in a right neighbourhood of $0$ ($h(\rho) \approx 0 $)  
than in a left neightbourhood of $1$ ($h(\rho) \approx 1 $), whereas assuming only that $h$ be bounded from below allows for many other instances like, e.g., a specular  behavior of $h$ around the extremal points of the domain~of~$f$. 

Returning now to \Ipblnew, we set
$$ g(\rho) : = h(\rho)  - \frac{\eps}2 \geq 0  \quad \hbox{ for all }\, \rho \in D(f_1),$$
and we reformulate our initial and boundary value problem as follows: 

\noindent to find $\mu$, $\rho$, and $\xi$, so as to solve
\Bsist
  \hskip-1cm& \bigl( \eps + 2g(\rho) \bigr) \, \dt\mu
  + \mu \, g'(\rho) \, \dt\rho
  - \div \bigl( \kamu\nabla\mu \bigr) = 0
  & \quad \hbox{in $\Omega \times (0,T)$,}
  \label{Iprima11}
  \\
   \hskip-1cm&\delta \dt\rho - \Delta\rho + \xi + f_2'(\rho) = \mu \,g'(\rho)
  \, , \ \hbox{ with } \
  \xi \in \beta (\rho),
   & \quad \hbox{in $\Omega \times (0,T)$,}
   \label{Iseconda11}
  \\
  \hskip-1cm& (\kamu\nabla\mu)\cdot\nu|_\Gamma = \dn\rho|_\Gamma = 0
  & \quad \hbox{on $\Gamma \times (0,T)$,}
  \label{Ibc11}
  \\
  \hskip-1cm& \mu(\cpto,0) = \muz
  \aand
  \rho(\cpto,0) = \rhoz
  & \quad \hbox{in $\Omega$.}
  \label{Icauchy11}
\Esist
\Accorpa\Ipblnewn Iprima11  Icauchy11
The \emph{global-in-time existence result} we derive is more general than in \cite{CGPS3}, for two reasons: because it holds when the potential $f$ includes a multivalued graph $\beta $ (possibly with vertical segments, e.g., for $f_1$ as in \eqref{ex2}) and only exploits the monotonicity property of $\beta$; and because the atom mobility $\kamu$ is allowed to depend in a generic nonlinear way on the chemical potential. {\pier Note that problem \Ipblnewn\ may become \emph{singular} with respect to $\rho$ due to the possible occurrence of singularities of $\beta$; on the other hand, it may be also \emph{degenerate} with respect to $\mu$ since $\kappa (\mu) $ is allowed to vanish at $\mu=0$.}

Our \emph{uniqueness result} is also more general than in \cite{CGPS3}, because nonsmooth potentials could not be handled with the technique there used, consisting in testing the difference of two equations \eqref{Iprima-3} by the time derivative of the difference of the two $\rho$ components; however, just as in \cite{CGPS3}, the proof is achieved under the assumption that $\kamu=$ a constant.
%

Some directions for future research have already been explored by us under less general circumstances than those considered here: 
the long time behavior of system \Ipbl\ and the structure of the relative 
omega-limit set have been analysed by us in \cite{CGPS3} and in \cite{CGPS4},
where we also dealt with the asymptotics of \Ipbl\ as $\eps \to 0$ and found a 
weaker solution in the singular limit. Moreover, in  \cite{CGPS5} and
\cite{CGS1} we studied two optimal control problems for systems
similar to \Ipbl: a distributed control problem in \cite{CGPS5} and a boundary control problem in \cite{CGS1}.
Finally, in \cite{CGPS7} we developed an existence theory for problem \Ipblnew\ when atom mobility  is allowed to depend on both $\mu$ and $\rho$.  
 
This paper is \organiz ed as follows.
In the next section, as anticipated, we discuss the physical features of the phase segregation model we adopt. In Section 3, we state our assumptions and results with the necessary mathematical accuracy; since here we do not take up asymptotic procedures, without loss of generality we set $\eps=1$ in \eqref{Iprima11} and $\delta=1$ in \eqref{Iseconda11}.
The existence of solutions to problem \Ipblnewn\ is proved in Section~\futuro\ref{Existence},  their regularity properties in the successive section. Our last Section~\futuro\ref{Uniqueness} is devoted to the uniqueness proof.
%
%


\section{{\pier Short} reasoned history of our mathematical model}
\label{reasoned}
\setcounter{equation}{0}

%

The nonstandard phase-field model \Ipbl\
can be regarded as a variant of the classic Cahn-Hilliard system for diffusion-driven phase segregation by atom rearrangement:
\Beq\label{CH}
 \dt\rho - \kappa \Delta \mu =0 \ , \qquad  \mu= - \Delta\rho + 
f'(\rho)
 {.}  
 \Eeq
Apart for the harmless choice $\kappa =1$ for the mobility modulus in \eqref{Iprima-3}, one finds in \Ipbl\ two  awkward nonlinear terms involving time derivatives. Usually, equations \eqref{CH} are combined in order to obtain the well-known \emph{Cahn-Hilliard equation}:
\Beq\label{CHe}
\dt\rho = \kappa \Delta (- \Delta\rho + f'(\rho)).
\Eeq
%
%

%
\step  Fried and Gurtin's generalization of Cahn-Hilliard equation

In \cite{FG,Gurtin},  a broad generalization of \eqref{CHe} was achieved by proposing the following:
\begin{enumerate}
 \item[(i)]
to interpret the second of \eqref{CH} as a \textit{balance of microforces:}
\Beq\label{balance}
\div\csib+\pi+\gamma=0{,}
\Eeq
where the distance microforce per unit volume 
is split into an internal part $\pi$ and an 
external part $\gamma$, and the contact microforce per unit area of a surface oriented by its normal $\nb$ is measured by $\csib\cdot\nb$ in terms of the \emph{microstress} vector  $\csib$;\footnote{In \cite{Fremond}, the microforce balance is stated under form of a principle of virtual powers for microscopic motions.}

\item[(ii)] 
to regard the first of \eqref{CH} as a 
\textit{balance law for the order parameter}:
\Beq\label{balorpam}
\partial_t\rho = - \div \hb + \sigma{,}
\Eeq
where the pair $(\hb ,  \sigma)$ is the \textit{inflow} of $\rho$; 

\item[(iii)] 
to demand that the constitutive choices for $\pi,\csib, \hb$, and the \emph{free energy density} $\psi$, {be} consistent in the sense of Coleman and Noll~\cite{CN} with an \emph{ad hoc} version of the Second Law of Continuum Thermodynamics:
\begin{equation}\label{dissipation}
\partial_t\psi +(\pi-\mu)\partial_t\rho-\csib\cdot\nabla(\partial_t\rho)+\hb\cdot\nabla\mu\leq 0 ,
\end{equation}
that is, a postulated ``dissipation inequality that {\js accommodates} diffusion'' (cf. equation~(3.6) in~\cite{Gurtin}). 
\end{enumerate}
In \cite{Gurtin}, the following set of constitutive prescriptions was shown to be consistent with~(iii):
\Beq\label{costi} 
\left\{
\begin{array}{c}
\psi = \widehat\psi(\rho,\nabla\rho),  \displaystyle \quad 
\widehat\pi(\rho,\nabla\rho,\mu)=\mu- \partial_\rho \widehat\psi(\rho,\nabla\rho), \displaystyle \quad 
\widehat\csib(\rho,\nabla\rho)=\partial_{\nabla\rho} \widehat\psi(\rho,\nabla\rho) \displaystyle
\end{array}
\right\} . 
\Eeq 
Moreover, it was presumed that
 \Beq\label{acca}
\hb = - \Mb\nabla \mu , \quad \hbox{with } \ \Mb=\widehat\Mb(\rho,\nabla\rho,\mu, \nabla\mu),
\Eeq
with the tensor-valued \emph{mobility mapping} $\widehat\Mb$ 
satisfying the  \emph{residual dissipation inequality}
\Beq
\nabla \mu\cdot \widehat\Mb(\rho,\nabla\rho,\mu, \nabla\mu) \nabla\mu \geq 0 . \non
\Eeq
With the help of \eqref{balance}, \eqref{balorpam}, \eqref{costi}, and $\eqref{acca}_1$,  a general equation for diffusive phase segregation processes is arrived at, namely,
%
\Beq\label{accag}
\dt\rho = \div\left(\Mb\nabla\left(\partial_\rho \widehat\psi(\rho,\nabla\rho)-\div\big(\partial_{\nabla\rho} \widehat\psi(\rho,\nabla\rho)\big)-\gamma\right)\right)+\sigma .
\Eeq
The classic Cahn-Hilliard equation \eqref{CHe} is obtained from \eqref{accag} by {taking} 
\Beq
\widehat\psi(\rho,\nabla\rho)= f(\rho)+\frac{1}{2}|\nabla\rho|^2,\qquad \Mb=\kappa \mathbf{1}, \label{constitutive}
\Eeq
and by letting the external distance microforce $\gamma$ and the order-parameter source term $\sigma$ be identically null.

\step An alternative generalization of Cahn-Hilliard equation

The Fried-Gurtin mod\-el was well accepted in the mathematical community. In 2006, a largely modified version of it was proposed \cite{Podio}: while the crucial step (i) was retained, both the order-parameter balance \eqref{balorpam} and the dissipation inequality \eqref{dissipation} were dropped and replaced, respectively, by
 the \emph{microenergy balance}
\Beq\label{energy}
\partial_t\varepsilon=e+w,\quad e:=-\div{\overline \hb}+{\overline \sigma},\quad w:=-\pi\,\partial_t\rho+\csib\cdot\nabla(\partial_t\rho){,}
\Eeq
and the \emph{microentropy imbalance}
\Beq\label{entropy}
\partial_t\eta\geq -\div\hb+\sigma,\quad \hb:=\mu{\overline \hb},\quad \sigma:=\mu\,{\overline \sigma}.
\Eeq
A further key feature of this new approach to modeling phase segregation by atomic rearrangement is that the \emph{microentropy inflow} $(\hb,\sigma)$ is deemed proportional to the \emph{microenergy inflow} $({\overline \hb},{\overline\sigma})$ through the \emph{chemical potential} $\mu$, a {positive} field; consistently, the free energy is defined~to~be
\Beq\label{freeenergy}
\psi:=\varepsilon-\mu^{-1}\eta,
\Eeq
with chemical potential playing the same role as  \emph{coldness} in the deduction of the heat equation.\footnote{As much as absolute temperature is a macroscopic measure of microscopic 
\emph{agitation}, its inverse - the coldness - measures microscopic \emph{quiet}; likewise, as argued in 
\cite{Podio}, {the} chemical potential can be seen as a macroscopic measure of microscopic \emph{organization}.}  

Combining (\ref{energy})-(\ref{freeenergy}) yields
\Beq\label{reduced}
\partial_t\psi\leq -\eta_{}\dt (\mu^{-1})+\mu^{-1}\,{\overline \hb}\cdot\nabla\mu-\pi\,\partial_t\rho+\csib\cdot\nabla(\partial_t\rho),
\Eeq
an inequality that replaces (\ref{dissipation}) in restricting  \emph{\`a la} Coleman-Noll {the possible} constitutive choices. 
On taking all of the constitutive mappings delivering $\pi,\csib,\eta$, and ${\overline \hb}$, 
depend{ent in principle} on  $\rho,\nabla\rho,\mu,\nabla\mu$, and on choosing
\Beq\label{constitutives}
\psi=\widehat\psi(\rho,\nabla\rho,\mu)=-\mu\,\rho+f(\rho)+\frac{1}{2}|\nabla\rho|^2,
\Eeq
compatibility with (\ref{reduced}) implies that we must have:
\Beq \label{cn}
\left\{
\begin{array}{c}
\widehat\pi(\rho,\nabla\rho,\mu)=-{\partial_{{\rho}} \widehat\psi(\rho,\nabla\rho,\mu)}~\displaystyle{=\mu-f'(\rho),}\\[0.2cm]
\widehat\csib(\rho,\nabla\rho,\mu)={\partial_{{\nabla\rho}} \widehat\psi(\rho,\nabla\rho,\mu)}=\nabla\rho, \displaystyle\\[0.2cm]
\widehat\eta(\rho,\nabla\rho,\mu)=\mu^2 \partial_{{\mu}} \widehat\psi(\rho,\nabla\rho,\mu)\displaystyle{=-\mu^2\rho}  
\end{array}
\right\}
\Eeq
together with
\Beq
\widehat{\overline \hb}(\rho,\nabla\rho,\mu,\nabla\mu) =  \widehat\Hb(\rho,\nabla\rho,\mu, \nabla\mu)\nabla \mu , \quad
\nabla \mu\cdot \widehat\Hb(\rho,\nabla\rho,\mu, \nabla\mu) \nabla\mu \geq 0 . \non
\Eeq
If {we now choose} for $\widehat\Hb$ the simplest expression $\Hb=\kappa \mathbf{1} $, implying a constant {and isotropic} mobility, and {if we once again} assume that the external distance microforce $\gamma$ and the source $\overline \sigma$ are null, {then}, with the use of (\ref{cn}) and \eqref{freeenergy}, the microforce balance (\ref{balance}) and the energy balance (\ref{energy}) become, respectively, 
\Beq\label{a}
 \Delta\rho+\mu- f'(\rho)  = 0
\Eeq
and
\Beq
 2\rho \,\dt\mu + \mu \, \dt\rho - \kappa \Delta \mu  = 0,   \label{secondanew}
\Eeq
a nonlinear system for the unknowns $\rho$ and $\mu$.

\step Insertion {\pier of the parameters} $\boldsymbol{\varepsilon}$ and 
$\boldsymbol \delta$

Let us compare systems \eqref{a}--\eqref{secondanew}
and \eqref{CH}.  {\pier Note that} \eqref{a} and $\eqref{CH}_2$ imply the same `static' relation between $\mu $ and $\rho$; instead,\eqref{secondanew}
is rather different from $\eqref{CH}_1$, for more than one reason: it is nonlinear; it features both time derivatives of $\rho$ and $\mu$; and,  in front of both $\dt\mu $ and 
$\dt\rho$ there are nonconstant factors that should remain nonnegative during the evolution. 
%
Thus, the system \eqref{a}--\eqref{secondanew} deserves a careful 
analysis. 

We begun by attacking the problem as it was, 
prompted to optimism by the successful outcome of a previous joint 
research effort~\cite{CGPS1, CGPS2} in which we tackled the system of Allen-Cahn 
type derived via the approach in \cite{Podio} 
for no-diffusion phase-segregation processes.
Unfortunately, {the evolution problem ruled by} \accorpa{a}{secondanew} turned out 
to be too difficult for us. Therefore, we decided to study its regularized 
version \Ipbl, where equations   \eqref{Iprima-3} and \eqref{Iseconda-3} are obtained by introducing the extra terms $\eps\, \partial_t\mu $ and $\delta\,\partial_t\rho$ in \eqref{secondanew} and \eqref{a}, respectively.
Of course, the positive coefficients $\eps$ and $\delta$ were intended 
to be made smaller and smaller by way of an asymptotic procedure to be set up after the solvability of the regularized system were proved.

 Mathematically, the introduction of the $\eps-$term is motivated by the 
desire to have a strictly positive coefficient as a factor 
of $\partial_t \mu$ in \eqref{secondanew},
so as to guarantee the parabolic structure of equation \eqref{Iprima-3}; on the other hand, the introduction of the  $\delta-$term transforms \eqref{a}
into an Allen-Cahn equation with source $\mu$, and is strongly reminiscent of a sort of regularization already employed in various approaches to the  
so-called \emph{viscous Cahn-Hilliard equations} (examples can be found in  
\cite{bcdgss, bds, GPS, mr, Ros} and in the references therein). 

It is also possible to make clear what additional physics the  
regularizing perturbations we introduced incorporate into the model. As to the term $\varepsilon\,\dt\mu$, it can be made to appear in the microenergy balance \eqref{Iprima-3} by modifying as follows the choice for the free energy in \eqref{constitutives}:
\begin{equation}
\label{fe-1}
\psi=-\mu\Big(\rho+\frac{\eps}{2}\Big)+f(\rho)+\frac{1}{2}|\nabla\rho|^2.
\end{equation}
As to the term $\delta\,\dt\rho$,   it is enough to note that all is needed to make that term appear in the microforce balance \eqref{Iseconda-3}
is to add $\dt\rho$ to the list of state variables we considered to analyze the constitutive consequences of \eqref{reduced}. This measure brings in the 
dissipation mechanism typical of Allen-Cahn nondiffusional segregation processes, where dissipation depends essentially on $(\dt\rho)^2$, in addition to Cahn-Hilliard's  $|\nabla\mu|^2-$ dissipation (cf.~\cite{Podio}); and it opens the way to \emph{splitting the distance microforce additively} into an equilibrium and a 
nonequilibrium part, with $\pi^{eq}=-
\partial_{{\rho}} \widehat\psi(\rho,\nabla\rho,\mu)=\mu-f'(\rho)$ the equilibrium part, just as in $\eqref{cn}_1$, and with $\pi^{neq}=-\delta\,\dt\rho$ the nonequilibrium part.

\section{Main results}
\label{MainResults}
\setcounter{equation}{0}

In this section, we state precisely the mathematical problem under 
investigation, fix our assumptions, and present our results.
Let $\Omega$ to be a bounded connected open set in $\erre^3$
with smooth boundary~$\Gamma$ {\js (the lower-dimensional cases can be treated
with minor changes)}. We also introduce a final time $T\in(0,+\infty)$
and set $Q:=\Omega\times(0,T)$. Moreover, we set for convenience:
\Beq
  V := \Huno,
  \quad H := \Ldue ,
  \aand
  W := \graffe{v\in\Hdue:\ \dn v = 0 \ \hbox{on $\Gamma$}},
  \label{defspazi}
\Eeq
and we endow these spaces with their standard norms,
for which we use a self-explanato\-ry notation 
like $\normaV\cpto$.
For $p\in[1,+\infty]$, we write $\norma\cpto_p$ for the usual norm
in~$L^p(\Omega)$. As no confusion can arise, the symbol $\norma\cpto_p$
is used for the norm in $L^p(Q)$ as well.
Moreover, any of the above symbols for the norms 
{\js  is used even for any power of these spaces}.
We remark that the embeddings $W\subset V\subset H$ are compact,
since $\Omega$ is bounded and smooth.
As $V$ is dense in~$H$, we can identify $H$ with
a subspace of $\Vp$ in the usual~way
(i.e.,~so as to have $\,_{\Vp}\<u,v>{}_V=(u,v)_H\,$
for every $u\in H$ and $v\in V$); the embedding $H\subset\Vp$ is also compact.

We are now concerned with the structural assumptions to set on our system.
As the chemical potential is {\js expected} to be at least nonnegative,
we assume that the function $\kappa$ is defined just 
for nonnegative {\js arguments}.
However, one could study the more general mathematical problem of 
finding solutions
whose component $\mu$ might change its sign.
In such a case, $\kappa$~has to be defined on the whole of $\erre$
and {\js must} satisfy similar assumptions.
We require~that:
\Bsist
  \hskip-1cm&& \hbox{$\kappa:[0,+\infty)\to\erre$ is continuous},
  \label{hpk}
  \\[1mm]
  \hskip-1cm&& \kmin,\kmax \in (0,+\infty)
  \aand
  \rmin \in [0,+\infty),
  \label{hpcost}
  \\[1mm]
  \hskip-1cm&& \kappa(r) \leq \kmax
  \quad \hbox{for every $r\geq0$}
  \aand
  \kappa(r) \geq \kmin
  \quad \hbox{for every $r\geq\rmin$},
  \label{hpkbis}
  \\[1mm]
  \hskip-1cm&& K(r) := \textstyle\int_0^r \kappa(s) \ds
  \quad \hbox{for $r\geq0$}; \quad
  \hbox{$K$ is strictly increasing},
  \label{hpK}
  \\[1mm]
  \hskip-1cm&& 
  {f = f_1 + f_2 \,, \quad f_1:\erre \to [0,+\infty], \quad f_2 :\erre \to \erre,} \quad   g:\erre \to [0,+\infty){,}
  \label{hpfg}
  \\[1mm]
  \hskip-1cm&& \hbox{$f_1$ is convex, proper, l.s.c., and $f_2$ and $g$ are $C^2$ functions}, 
  \label{hpfuno}
  \\[1mm]
  \hskip-1cm&& \hbox{$f_2'$, $g$, and $g'$ are Lipschitz continuous},
  \label{hpfdueg}
  \\[1mm]
  \hskip-1cm&& \beta := \partial f_1 \aand
  \pi := f_2' \,.
  \label{defbp}
\Esist
\Accorpa\Hpstruttura hpk defbp
In the {\js following}, $D(f_1)$ and $D(\beta)\, {\pier (\subseteq D(f_1))}$ 
denote the effective domains
of~$f_1$ and~$\beta$, respectively.

\Brem
\label{Remstruttura}
We observe that our assumptions on $f$ and $\kappa$
allow {\js for} strong {\pier singularities (at the boundary of $D(\beta)$)} 
and {\js a} possible degeneracy {\pier (in a right neighbourhood of $0$)} in the equations for 
$\rho$ and~$\mu$, respectively.
The former fact is clear.
As far as the latter is concerned,
we note that \eqref{hpK} is satisfied if and only if 
$\kappa$~is nonnegative and the set
where $\kappa$ vanishes {\js has empty interior}.
So, $\rmin=0$ means uniform parabolicity for equation~\eqref{Iprima}.
On the contrary, if $\rmin>0$, the equation can degenerate,
e.g., at the {\js origin} (or~even in rather big set of small values).
An example is given by $\kappa(r)=\tanh r^{m-1}$ with $m>1$.
In {\js this} case, \eqref{Iprima} roughly behaves 
like the porous medium equation (slow diffusion) 
in the region where $\mu$ is small.
\Erem

{\pier \Brem
\label{justif}
It is known that any proper, convex and \lsc\ function is bounded from below by an 
affine function (see, e.g., \cite[Prop.~2.1, p.~51]{Barbu}). Hence,
our  assumption $f_1 \geq 0 $ looks reasonable, because one can suitably 
modify the smooth perturbation~$f_2$ by adding a straight line. 
{\js Moreover,} the other positivity condition{,} $g\geq 0${,} 
is just needed on the set~$D(\beta)$, while $g$ can {\pier take 
negative} values outside of $D(\beta)$. {\js Finally, \eqref{hpfdueg}
implies that{\pier ,} within the range of relevant values of $r$, the functions 
{\pier $f_2'(r)$, $g(r)$, and $g'(r)$} grow at most linearly with respect to $|r|$, 
while {\pier $f_2(r)$} grows at most quadratically in~$|r|$.}   
\Erem
}

{\js For the initial data, we postulate:}
\Bsist
  && \muz \in V , \quad
  \rhoz \in W , \quad
  \muz \geq 0
  \aand \rhoz \in D(\beta) \quad \aeO ,
  \label{hpzero}
  \\
  && \hbox{and there exists some $\xi_0\in H$ such that $\xi_0\in\beta(\rhoz)$ \aeO}.
  \label{hprhozbis}
\Esist
\Accorpa\Hpdati hpzero hprhozbis
{\js Since} $f_1$ is convex and $f_2$ is smooth,
the above assumptions entail $f(\rhoz)\in\Luno$.

Now, we introduce the a priori regularity {\js that we require 
from any solution $(\mu,\rho,\xi)$ to our problem}.
Note that {\js equation~\eqref{Iseconda}
reduces for any given~$\mu$ to} a rather standard {\pier phase-field} equation.
Therefore, it is natural to look for pairs $(\rho,\xi)$ that~satisfy
\Bsist
  &\rho \in \W{1,\infty}H \cap \H1V \cap \L\infty W ,&
  \label{regrho}
  \\[1mm]
  & \xi \in \L\infty H ,&
  \label{regxi}
\Esist
and {\js to} solve the subproblem in a strong form, namely
\Accorpa\Regrhoxi regrho regxi
\Bsist
  & \dt\rho - \Delta\rho + \xi + \pi(\rho) = \mu \, g'(\rho)
  \aand \xi \in \beta(\rho)
  & \quad \aeQ,
  \label{seconda}
  \\
  & 
  \rho(0) = \rhoz
  & \quad \aeO .
  \label{cauchyrho}
\Esist
\Accorpa\Pblrhoxi seconda cauchyrho
We note that \eqref{regrho} also contains the Neumann boundary condition 
{\pier for $\rho$  (see~\eqref{defspazi} for the definition of $W$)}.
On the contrary, the situation is different for the component~$\mu$.
In the case of uniform parabolicity, i.e., {\js if} $\rmin=0$,
the coefficient $\kamu$ is bounded away from zero,
and we can require~that
\Bsist
  &\mu \in \H1H \cap \L\infty V, \quad
  \mu \geq 0 \quad \aeQ , &
  \label{regmu}
  \\
  & \div \bigl( \kamu\nabla\mu \bigr) \in \L2H ,&
  \label{regdiv}
\Esist
\Accorpa\Regsoluz regrho regdiv
and that $\mu$ satisfy 
\Bsist
  \hskip-1.5cm& \displaystyle \iO \bigl( \coefft \bigr) \, \dt\mu(t) \, v
  + \iO \mu(t) \, g'(\rho(t)) \, \dt\rho(t) \, v
  + \iO \kappa(\mu(t)) \nabla\mu(t) \cdot \nabla v = 0 &
  \non
  \\
   \hskip-1.5cm& \hskip6cm\hbox{for every $v\in V$ and \aat},&  
  \label{prima}
  \\
   \hskip-1.5cm& \mu(0) = \muz
  \quad \aeO .&  \vphantom \int
  \label{cauchymu}
\Esist
\Accorpa\Pblmu prima cauchymu
Thus, the equation holds in a strong sense, i.e.,
\Beq
  \bigl( \coeff \bigr) \, \dt\mu
  + \mu \, g'(\rho) \, \dt\rho
  - \div \bigl( \kamu\nabla\mu \bigr) = 0
  \quad \aeQ , 
  \label{primaaeQ}
\Eeq
while the Neumann boundary condition is understood in the usual weak sense.
Furthermore, we observe that \accorpa{regmu}{prima} imply {\pier further regularity 
for $\mu$} whenever $\kappa$ is smoother, thanks to
the regularity theory of {\pier quasilinear} elliptic equations.

If instead we allow $\rmin$ to be positive,
such a formulation is too strong,
since no sufficient information on the gradient $\nabla\mu$
can be obtained.
As a consequence, the same happens for the time derivative~$\dt\mu$.
For that reason, we rewrite equation~\eqref{primaaeQ} in the different form
\Beq
  \dt \bigl( \coeff \mu \bigr) - \mu g'(\mu) \dt\rho - \Delta K(\mu) = 0.
  \label{primadistr}
\Eeq
More precisely, we also account for the initial and Neumann boundary conditions
and accordingly rewrite~\Pblmu.
In conclusion, we require the lower regularity
\Bsist
  &\mu \in \L\infty H , \quad
  \mu \geq 0 \quad \aeQ , \quad
  K(\mu) \in \H1H \cap \L\infty V,  &
  \qquad 
  \label{wregmua}
  \\[1mm]
  & (\coeff) \mu \in \H1\Vp, &
  \label{wregmub}
\Esist
\Accorpa\Wregmu wregmua wregmub
and replace \Pblmu~by
\Bsist
  & \displaystyle \< \dt \bigl( (\coeff)\mu \bigr)(t) , v >
  - \iO \bigl( \mu g'(\rho)\dt\rho \bigr)(t) \, v
  + \iO \nabla K(\mu(t)) \cdot \nabla v
  = 0 &
  \non
  \\
  & \hskip5cm\hbox{for every $v\in V$ and \aat},&
  \label{wprima}
  \\
  & \bigl( (\coeff)\mu \bigr)(0) = \bigl( 1 + 2g(\rhoz ) \bigr) \, \muz.&
  \vphantom \int
  \label{wcauchymu}
\Esist
\Accorpa\Wpblmu wprima wcauchymu
In {\js this situation}, \eqref{primadistr} is satisfied in the sense of distributions, only.

\Brem
\label{Middleterm}
We observe that even the middle term of \eqref{wprima} is meaningful,
as we immediately see.
First, we note~that
\Beq
  \rho \in \C0{\Cx0} = \CQ,
  \label{rhocont}
\Eeq
directly from \eqref{regrho} and the compact embedding $W\subset\Cx0$
(see, e.g., \cite[Sect.~8, Cor.~4]{Simon}),
whence $g'(\rho)\in \CQ$.
Next, \eqref{wregmua} and the embedding $V\subset\Lq$
imply that $K(\mu)\in\L\infty\Lq$, whence also $\mu\in\L\infty\Lq$,
since $K(r)$ behaves like $r$ for big~$|r|$ by \eqref{hpkbis}.
Finally, $\dt\rho\in\L\infty H$.
Therefore, $\mu g'(\rho)\dt\rho\in\L\infty{\Lx{4/3}}$.
On the other hand, $v\in\Lq$ whenever $v\in V$.
\Erem

\Brem
\label{Weakcont}
Note that \eqref{wcauchymu} makes sense 
because $(\coeff)\mu\in\C0\Vp$ (due to~\eqref{wregmub}).
However, by accounting for \eqref{regrho} and the regularity of~$g$,
we see that \eqref{wcauchymu} can be read in the simpler form \eqref{cauchymu}
also in this case.
Indeed, the function $(\coeff)\mu$ also belongs to $\L\infty H$.
As is well known (and easy to prove),
this implies that it actually is an $H$-valued function
which is weakly continuous, in addition.
It easily follows that $\mu$ enjoys the same property.
\Erem

Here are our results.
The first {\js two} state that the strong and weak formulations
are equivalent in the case $\rmin=0$
and that there exists a weak solution in the general case.
Due to the former, the latter also
proves the existence of a strong solution if $\rmin=0$.
{\pier Both results {\js will be} proved in Section~\ref{Existence}.}

\Bprop
\label{Equivalenza}
Assume \Hpstruttura, \Hpdati, and $\rmin=0$.
Then, any triplet $(\mu,\rho,\xi)$ satisfing \Regrhoxi, \Wregmu\ and solving
problem~\Pblrhoxi{\pier , \Wpblmu}\
also satisfies \accorpa{regmu}{cauchymu}.
\Eprop

\Bthm
\label{Esistenza}
Assume \Hpstruttura\ and \Hpdati.
Then, there exists at least {\js one} triplet $(\mu,\rho,\xi)$ satisfing \Regrhoxi, \Wregmu\ and solving
problem~\Pblrhoxi{\pier , \Wpblmu}.
\Ethm

We notice that no further assumptions are needed
to ensure boundedness for~$\rho$, due to~\eqref{rhocont}.
As far as the first component is concerned, we have the following boundedness result.

\Bthm
\label{Mubdd}
Assume \Hpstruttura, \Hpdati, and let~${\pier \muz\in L^\infty (\Omega)}$. 
Then, the component $\mu$ of any triplet $(\mu,\rho,\xi)$ 
satisfing \Regrhoxi, \Wregmu\ and solving
problem~\Pblrhoxi\ and \Wpblmu\ is {\pier essentially} bounded.
\Ethm

The next results hold if we assume
that $\kappa$ is constant.
We notice that we could weaken {\js this} assumption 
in our regularity result,
while we are not able to do the same
as far as uniqueness is concerned, unfortunately.
In order to simplify the regularity proof,
we take $\kappa=1$ at once.
In the forthcoming Remark~\futuro\ref{Ancorareg},
we {\js will sketch how to deduce even further} regularity.

\Bthm
\label{Piuregmu}
Assume \Hpstruttura, \Hpdati, $\muz\in W$, and $\kappa=1$.
Then, any triplet $(\mu,\rho,\xi)$ satisfing \Regrhoxi, \eqref{regmu} 
and solving problem~\Pblrhoxi\ and \Pblmu\
enjoys the regularity property
\Beq
  \mu \in \W{1,p}H \cap \L pW
  \quad \hbox{for every $p\in[1,+\infty)$}.
  \label{piuregmu}
\Eeq
\Ethm

\Bthm
\label{Unicita}
Assume \Hpstruttura, \Hpdati, $\muz\in W$, and $\kappa=1$.
Then, the triplet $(\mu,\rho,\xi)$ satisfing \Regrhoxi, \eqref{regmu} 
and solving problem~\Pblrhoxi\ and \Pblmu\ is unique.
\Ethm

Throughout the paper,
we account for the \wk\ embedding $V\subset\Lx p$
for $1\leq p\leq 6$
and the related Sobolev inequality:
\Beq
  \norma v_p \leq C \normaV v
  \quad \hbox{for every $v\in V$ and $1\leq p \leq 6$,}
  \label{sobolev}
\Eeq
where $C$ depends on~$\Omega$ only.
Moreover, we recall that
the embeddings $V\subset\Lq$ 
(more generally $V\subset\Lx p$ with $p<6$) 
and $W\subset\Cx0$ are compact
and use the corresponding inequality
\Beq
  \norma v_4 \leq \eps \normaH{\nabla v} + C_\eps \normaH v
  \quad \hbox{for every $v\in V$ and $\eps>0$,}
  \label{compact}
\Eeq
where $C_\eps$ depends on $\Omega$ and~$\eps$, only.
Furthermore, we {\js make repeated} use of the notation
\Beq
  Q_t := \Omega \times (0,t)
  \quad \hbox{for $t\in[0,T],$}
  \label{defQt}
\Eeq
and of the \wk\ \holder\ inequality
and the elementary Young inequality
\Beq
  ab \leq \eps a^2 + \frac 1{4\eps} \, b^2
  \quad \hbox{for every $a,b\geq 0$ and $\eps>0$}.
  \label{young}
\Eeq
Finally, throughout the paper
we use a small-case italic $c$ for different constants that
may only depend 
on~$\Omega$, the final time~$T$, the shape of the nonlinearities $f$ and~$g$, 
and the properties of the data involved in the statements at hand; 
a~notation like~$c_\eps$ signals a constant that depends also on the parameter~$\eps$. 
The reader should keep in mind that the meaning of $c$ and $c_\eps$ might
change from line to line and even in the same chain of inequalities, 
whereas those constants {\js that} we need to refer to are always denoted by 
capital letters, just like $C$ in~\eqref{sobolev}.


\section{Existence}
\label{Existence}
\setcounter{equation}{0}

In this section, we {\pier first show} the equivalence result stated in Proposition~\ref{Equivalenza}{\pier . Then, we prove}
Theorem~\ref{Esistenza}, which ensures the existence of a weak solution.

\proofstep
Proof of Proposition~\ref{Equivalenza}

As $\rmin=0$, we have $\kappa(r)\geq\kmin$ for every $r\geq0$.
This implies that the inverse function $K^{-1}:[0,+\infty)\to[0,+\infty)$
is \Lip\ continuous.
Hence, \eqref{wregmua} implies that
\Beq
  \mu = K^{-1}(K(\mu)) \in \H1H \cap \L\infty V,
  \non
\Eeq
i.e., that \eqref{regmu} holds.
In particular, we can write
\Beq
  \nabla K(\mu) = \kamu \nabla\mu
  \aand
  \dt \bigl( (\coeff)\mu \bigr)
  = \mu \dt (\coeff) + (\coeff) \dt\mu
  \non
\Eeq
and thus replace the weak formulation by the strong one.
Next, we note that \eqref{prima} implies that
\eqref{primaaeQ} holds in the sense of distribution,
whence \eqref{regdiv} follows by comparison.
Finally, \eqref{cauchymu} holds even in the general case, 
as we have observed in Remark~\ref{Weakcont}.\QED

\medskip

Now we prove Theorem~\ref{Esistenza}.
Even though our proof closely follows the argument{\pier ation} of~\cite{CGPS3},
we present the whole {\pier procedure} and sometimes give some detail,
since the changes with respect to the quoted paper are spread in the whole calculation.
The starting point is an approximating problem
which is still based on a time delay in the \rhs\
of~\eqref{seconda}.
Namely, we define the translation operator $\T:\L1H\to\L1H$
depending
on a time step $\tau>0$ by setting, for $v\in\L1H$ and \aat,
\Beq
  (\T v)(t) := v(t-\tau)
  \quad \hbox{if $t>\tau$}
  \aand
  (\T v)(t) := \muz
  \quad \hbox{if $t<\tau$}
  \label{defT}
\Eeq 
(but the same notation $\T v$ will be {\js used even for functions $v$
that are} defined in some subinterval $[0,T']$ of~$[0,T]$),
and replace $\mu$ by $\T\mu$ in~\eqref{seconda}, essentially.
However, we modify the equation for $\mu$ at the same time.
Precisely, we force uniform parabolicity 
and allow the solution
to take negative values, if possible.
To do that, we define $\kappat:\erre\to\erre$ 
and the related function $\Kt$ used later on~by
\Beq
  \kappat(r) := \kappa(|r|)+\tau
  \aand
  \Kt(r) := \int_0^r \kappat(s) \ds
  \quad \hbox{for $r\in\erre$}.
  \label{defkt}
\Eeq
So, the approximating problem consists of the following equations
\Bsist
  & \bigl( \coefftau \bigr) \, \dt\mut
  + \mut \, g'(\rhot) \, \dt\rhot
  - \div \bigl( \kamut\nabla\mut \bigr) = 0
  & \quad \aeQ,
  \label{primatau}
  \\[1mm]
  & \dt\rhot - \Delta\rhot + \xit + \pi(\rhot)
  = (\T\mut) \, g'(\rhot)
  \aand \xit \in \beta(\rhot)
  & \quad \aeQ,
  \label{secondatau}
\Esist
complemented with the homogeneous Neumann boundary conditions
for both $\mut$ and $\rhot$
and the initial conditions $\mut(0)=\muz$ and $\rhot(0)=\rhoz$.
For convenience, we allow $\tau$ to take just discrete values,
namely, $\tau=T/N$, where $N$ is any positive integer.\QED

\Blem
\label{EsistPbltau}
The approximating problem has a solution $(\mut,\rhot,\xit)$ satisfying
the {\js analogues} of \Regrhoxi\ and \accorpa{regmu}{regdiv}.
\Elem

\Bdim
We {\pier just} give a sketch.
As in~\cite{CGPS3}, we inductively solve $N$ problems
on the time intervals $I_n=[0,t_n]$, $n=1,\dots,N$,
by constructing the solution directly on the whole of~$I_n$ at each step.
Namely, given $\munmu$, which is defined in $\Omega\times I_{n-1}$,
we note that $\T\munmu$ is well defined and known 
in $\Omega\times I_n$ (even in the starting case $n=1$) 
and solves the boundary value problem for $\rhon$
given by the {\pier phase-field} equations
\Beq
  \dt\rhon - \Delta\rhon + \xin + \pi(\rhon)
  = (\T\mun) \, g'(\rhon)
  \aand \xin \in \beta(\rhon)
  \quad \hbox{in $\Omega\times I_n$},
  \label{secondan}
\Eeq
complemented with the boundary and initial conditions just mentioned for~$\rhot$.
Such a problem is quite standard and has a unique solution 
in a proper (rather weak) functional {\js analytic} framework.
Once $\rhon$ is constructed, we solve the parabolic equation
\Beq
  \bigl( \coeffn \bigr) \, \dt\mun
  + \mun \, g'(\rhon) \, \dt\rhon
  - \div \bigl( \kappat(\mun)\nabla\mun \bigr) = 0
  \quad \hbox{in $\Omega\times I_n$},
  \label{priman}
\Eeq
{\js together} with the boundary and initial conditions prescribed for~$\mut$.
We note that $g\geq0$ and $\kappat(r)\geq\tau$ for every $r\in\erre$,
so that the equation is uniformly parabolic.
Therefore, the problem to be solved has a unique solution in a proper space
provided that the coefficient $g'(\rhon)\dt\rhon$ is not too irregular.
So, we should prove that, step by step, we get the right regularity
for $\rhon$ and~$\mun$.
This could be done by induction, as in~\cite{CGPS3}, with some modifications
due to our more general framework.
We omit this detail and just observe that the needed a~priori estimates 
are close (and even simpler, since $\tau$ is fixed here)
to~the ones {\js performed} later on in order to let $\tau$ {\js tend} to zero.
The final point is $\mun\geq0$.
We give the proof in detail.
We multiply equation~\eqref{priman} by $-\mun^-:=-(-\mun)^+$,
the negative part of~$\mun$,
and integrate over~$Q_t$ with any $t\in I_n$.
We observe that
\Beq
  \bigl[ \bigl( \coeffnt \bigr) \, \dt\mun + \mun \, g(\rhon) \, \dt\rhon \bigr] \, (-\mun^-)
  = \frac 12 \, \dt \bigl( (\coeffn) \, |\mun^-|^2 \bigr).
  \non
\Eeq
Hence, by using $\muz\geq0$, 
and owing to the boundary condition, we have
\Bsist
  \frac 12 \iO (\coeffnt) \, |\mun^-(t)|^2
  + \intQt \kappat(\mun) |\nabla\mun^-|^2 
  = 0 
  \quad \hbox{for every $t\in I_n$}.
  \non
\Esist
Since $g$ and $\kappat$ are nonnegative,
this implies $\mun^-=0$, that is, $\mun\geq0$ a.e.\ in~$\Omega\times I_n$.
Once all this is checked, the finite sequence
$(\mun,\rhon,\xin)$, $n=1,\dots,N$, {\js is actually} constructed, 
and it is clear that 
a solution to the approximating problem we are looking for
is {\js  obained by simply }taking $n=N$.
\Edim

{\pier Although} the solution to the approximating problem is unique,
we do not need uniqueness in the {\js following}
and just {\pier fix} a solution $(\mut,\rhot,\xit)$ for each~$\tau$.
Our aim is to let $\tau$ {\js tend} to zero in order to obtain
a solution as stated in Theorem~\ref{Esistenza}.
Our proof uses compactness arguments and thus relies on a number of 
uniform (with respect to~$\tau$) a~priori estimates.
Clearly, in performing them, we can take $\tau$ as small as we desire,
and it will be suitable to assume that $\tau\leq\kmax$.
In order to make the formulas more readable, 
we shall omit the index~$\tau$ in the calculations, 
waiting for writing $\mut$ and $\rhot$ only when each estimate is established.

\step 
First a priori estimate

{\pier Let us test \eqref{primatau} by $\mut$ and point out} that
\Beq
  \bigl[ \bigl( \coeff \bigr) \, \dt\mu
  + \mu \, g'(\rho) \, \dt\rho \bigr] \mu
  = \frac 12 \, \dt \bigl[ (\coeff) \mu^2 \bigr].
  \non
\Eeq
{\pier Therefore}, by integrating over~$(0,t)$, where $t\in[0,T]$ is arbitrary,
we obtain
\Beq
  \iO \bigl( \coefft \bigr) |\mu(t)|^2 + \intQt \kappat(\mu) |\nabla\mu|^2
  = \iO (1+2g(\rhoz)) \muz^2 \,.
  \non
\Eeq
Hence, we recall that $g\geq0$ and observe that 
$\kappat^2(r)\leq 2\kmax\kappat(r)$ for every $r\in\erre$
by \eqref{hpkbis} and $\tau\leq\kmax$.
We conclude that
\Beq
  \norma\mut_{\L\infty H} + \norma{\Kt(\mut)}_{\L2V} \leq c .
  \label{primastima}
\Eeq
Actually, we have proved more, namely
\Beq
  \norma{K^*_\tau(\mut)}_{\L2V} \leq c 
  \quad \hbox{where} \quad
  K^*_\tau(r) := \int_0^r (\kappat(s))^{1/2} \ds
  \quad \hbox{for $r\in\erre$}.
  \non
\Eeq
Moreover, we observe that $K$ has a linear growth,
so that \eqref{primastima} also yields
\Beq
  \norma{\Kt(\mut)}_{\L\infty H} \leq c .
  \label{daprimaK}
\Eeq
{\pier An implication of \accorpa{primastima}{daprimaK}, along with \eqref{defT} and \eqref{hpzero}, is}
\Beq
 {\pier  \norma{\T \mut}_{\L\infty H} + \norma{\T \Kt(\mut)}_{\L\infty H\cap \L2V} \leq c .}
  \label{pristitau} 
\Eeq

\step Consequence

The Sobolev inequality~\eqref{sobolev} and estimate \eqref{primastima} imply that
\Beq
  \norma{\Kt(\mut)}_{\L2{\Lx6}} \leq c .
  \non
\Eeq
On the other hand, \eqref{hpkbis} implies that $\Kt(r)\geq\kmin r-c$ for every $r\geq0$.
We deduce that
\Beq
  \norma\mut_{\L2{\Lx6}} \leq c
  \label{daprima}
\Eeq

\step 
Second a priori estimate

{\pier Let us add $\rhot$ on both sides of \eqref{secondatau} and test by $\dt\rhot$.
We have that 
\Bsist
  && \intQt |\dt\rho|^2
  + {\pier \frac 12 \, \norma{\rho(t)}_V^2 }
  + \iO f_1(\rho(t))
  \non
  \\
  && = \frac 12 {\js \|\rhoz\|^2_V}
  + \iO f(\rhoz)
  + {\pier \frac 12 \iO \tonde{\rho^2 (t) {\js -} 2 f_2(\rho(t))} }
  + \intQt g'(\rho) (\T\mu) \dt\rho 
  \non
  \\
  && \leq c + {\pier c \iO |\rho (t)|^2} + \frac 14 \intQt |\dt\rho|^2
  + c \norma{\T\mu}_{\L\infty H}^2 ,
  \non
\Esist
for every $t\in[0,T]$. In view of the chain rule and Young's inequality \eqref{young}, we observe that
\Beq
  c \iO |\rho(t)|^2 
  \leq c \iO |\rhoz|^2 + \frac 14  \intQt |\dt\rho|^2  
  + c\int_0^t \norma{\rho(s)}_H^2\, ds.
  \non
\Eeq
Hence, as $f_1$ is nonnegative, on account of \eqref{pristitau}, and 
with the help of the Gronwall lemma, we deduce that
\Beq
  \intQt |\dt\rho|^2
  +  \norma{\rho(t)}_V^2
  + \iO  f_1(\rho(t)) \leq c.
\non
\Eeq
Therefore, we obtain:
\Beq
  \norma\rhot_{\H1H\cap\L\infty V} \leq c
  \aand
  \norma{f(\rhot)}_{\L\infty\Luno} \leq c.
  \label{secondastima}
\Eeq
}

\step 
Third a priori estimate

We rewrite \eqref{secondatau} as 
\Beq
  -\Delta\rho + \beta(\rho)
  \ni -\dt\rho - \pi(\rho) + (\T\mu)g'(\rho) 
  \non
\Eeq
and note that the \rhs\ is bounded in $\L2H$, 
thanks to the \Lip\ continuity of $\pi$ and $g'$ and 
{\pier to} the previous estimates.
By a standard argument
({\pier {\js formally} test} by either $-\Delta\rho$ or $\beta(\rho)$
and use  the regularity theory for elliptic equations),
{\pier we first recover {\js that} 
\Beq
  \normaH{\Delta\rho(s)}^2 + \normaH{\xi(s)}^2
  \leq  2 \normaH{-\dt\rho(s) - \pi(\rho(s)) +  ((\T\mu)g'(\rho))(s)}^2
  \label{pier1}
\Eeq
for a.a. $s\in (0,T),$ and~finally} conclude that
\Beq
  \norma\rhot_{\L2W} \leq c
  \aand
  \norma\xit_{\L2H} \leq c.
  \label{terzastima}
\Eeq

\step 
Fourth a priori estimate

Our aim is to improve {\js the} estimates~\eqref{secondastima} and~\eqref{terzastima}.
To {\js this end}, we proceed formally, at least at the beginning,
for the sake of simplicity.
However, {\js this} procedure could be made rigorous
by suitably \regulariz ing equation \eqref{secondatau}
(with respect to $\rho$, only, i.e., keeping $\mu$ fixed),
the main tool being the Yosida \regulariz ation of maximal monotone operators
(see, e.g., \cite[p.~28]{Brezis}; 
see also the proof of Lemma~3.1 of~\cite{CGPS3} for a further \regulariz ation).
Such a theory yields, in particular, the estimate
\Beq
  \normaH{\dt u(0)}
  \leq \normaH{\psi(0)+\Delta\rhoz}
  + \min_{\eta\in\beta(\rhoz)} \normaH\eta
  \label{perquartaZ}
\Eeq
for the unique solution $(u,\omega)$ to the equations
\Beq
  \dt u - \Delta u + \omega = \psi := g'(\rho) \T\mu - \pi(\rho)
  \aand
  \omega \in \beta(u) ,
  \non
\Eeq
complemented {\js with} the same initial and boundary conditions {\js as those} prescribed for~$\rho$.
{\pier Note that in \eqref{perquartaZ} $\beta$ is understood 
as the induced maximal 
monotone operator from $H$ to $H$. Observe} that $(u,\omega)=(\rho,\xi)${\pier ; then}  {\js the
application of \eqref{perquartaZ}, in combination} with our assumptions 
on~$\rhoz$ (see \eqref{hprhozbis}, in particular), {\pier leads to} 
\Beq
  \normaH{\dt\rhot(0)}
  \leq c \bigl( \normaH\muz + \normaW\rhoz + 1 + \normaH{\xi_0} \bigl)
  = c .
  \label{finequartaZ}
\Eeq
We use \eqref{finequartaZ} in the {\js subsequent calculation,
where we proceed formally, as announced
(however, our procedure becomes completely rigorous after a while).
In~particular, we write $\beta(\rho)$ in place of $\xi$ 
and treat $\beta$ as {\js if it were} a smooth function.
We differentiate \eqref{secondatau} with respect to time
and test the resulting} equation by $\dt\rho$.
We obtain{\pier :}
\Bsist
  && \frac 12 \iO |\dt\rho(t)|^2
  + \intQt |\nabla\dt\rho|^2
  + \intQt \beta'(\rho) |\dt\rho|^2
  \non
  \\
  && = \frac 12 \iO |(\dt\rho)(0)|^2 
  - \intQt \pi'(\rho) |\dt\rho|^2
  + \intQt g''(\rho) {\pier (\T\mu)} |\dt\rho|^2
  \qquad
  \non
  \\
  &&\quad {} 
  + \intQt g'(\rho) \dt(\T\mu) \, \dt\rho
  \non
  \\
  && \leq \frac 12 \iO |(\dt\rho)(0)|^2 
  + c \intQt (1+{\pier \T \mu}) |\dt\rho|^2
  + \intQt g'(\rho) \dt(\T\mu) \, \dt\rho .
  \label{perquarta}
\Esist
We treat each term on the {\rhs} separately.
The first one is estimated by~\eqref{finequartaZ}.
In order to deal with the second one,
we first account for the \holder\ inequality.
Then, we also {\js invoke} \eqref{primastima}, 
the compact embedding $V\subset\Lq$ (see \eqref{compact}),
and \eqref{secondastima}. {\pier For every~$\eps>0$ we~infer that}
\Bsist
  && \intQt (1+{\pier \T \mu}) |\dt\rho|^2
  \leq \iot \normaH{1+{\pier (\T \mu )}(s)} \norma{\dt\rho(s)}_4^2 \ds
  \non
  \\
  && \leq c \iot \norma{\dt\rho(s)}_4^2 \ds
  \leq \eps \intQt |\nabla\dt\rho|^2
  + c_\eps \intQt |\dt\rho|^2
  \qquad
  \non
  \\
  && \leq \eps \intQt |\nabla\dt\rho|^2
  + c_\eps \,. 
  \label{finequartaA}
\Esist
The estimate of the last term of~\eqref{perquarta} needs much more work.
We recall that $\T\mu$ is constant with respect to time on the interval~$(0,\tau)$
and first compute $\dt\mu$ from~\eqref{primatau}.
Then we integrate by parts and repeatedly {\pier exploit} 
the \holder, Sobolev, and Young inequalities.
We obtain{\pier :}
\Bsist
  \hskip-1cm && \intQt g'(\rho) \dt(\T\mu) \, \dt\rho
  = \intQtmt \dt\mu(s) \, g'(\rho({\pier s+\tau})) \dt\rho(s+\tau) \ds
  \non
  \\
  \hskip-1cm && = \intQtmt \frac 1{\coeffs} \,
  \bigl[ \div \bigl({\pier\kappat(\mu)(s)}\nabla\mu(s)\bigr) - \mu(s) g'(\rho(s)) \dt\rho(s) \bigr]
  \dt{\pier g(\rho(s+\tau))} \ds
  \non
  \\
  \hskip-1cm && = \intQtmt {\pier\kappat(\mu)(s)} \nabla\mu(s) \cdot \nabla 
  \frac {\dt {\pier g(\rho(s+\tau))}} {\coeffs} \ds
  \non
  \\
  \hskip-1cm && \qquad {}
  - \intQtmt \frac {g'(\rho(s)) \, {\pier g'(\rho(s+\tau))} }{\coeffs} \, \mu(s) \dt\rho(s) \dt\rho(s+\tau) \ds, 
  \label{perquartaB}
\Esist
and now treat the last two integrals separately,
by accounting for our structural assumptions.
We have {\pier
\Bsist
  && \intQtmt {\pier\kappat(\mu)(s)} \nabla\mu(s) \cdot \nabla \frac {\dt  g(\rho(s+\tau))}  {\coeffs} \ds \non
  \\
  &&= \intQtmt \nabla{\pier\Kt(\mu)(s)} \cdot\nabla \, 
  \frac{g'(\rho(s+\tau)) \dt \rho(s+\tau)}{\coeffs} \ds
  \non
  \\
  && \leq c \intQtmt |\nabla{\pier\Kt(\mu)(s)}| \, |\nabla\dt\rho(s+\tau)| \ds
  \non
  \\
  && \quad {}
  + c \intQtmt |\nabla{\pier\Kt(\mu)(s)}| \, |\nabla\rho(s)| \, |\dt\rho(s+\tau)| \ds
  \non
  \\
  && \quad {}
  + c \intQtmt |\nabla{\pier\Kt(\mu)(s)}| \, |\nabla\rho(s +\tau)| \, |\dt\rho(s+\tau)| \ds. \label{perquartaBA}
\Esist
}
{\pier For every $\eps\in(0,1)$ there holds}
\Bsist
  && \intQtmt |\nabla{\pier\Kt(\mu)(s)}| \, |\nabla\dt\rho(s+\tau)| \ds
  \non
  \\
  && 
  {}\leq \eps \intQt |\nabla\dt\rho|^2
  + c_\eps \intQt |\nabla\Kt(\mu)|^2 
  \, \leq \, \eps \intQt |\nabla\dt\rho|^2 + c_\eps,
  \label{perquartaBAA}
\Esist
thanks to~\eqref{primastima}.
On the other hand, {\pier we also have}
\Bsist
  && \intQtmt |\nabla{\pier\Kt(\mu)(s)}| \, |\nabla\rho(s)| \, |\dt\rho(s+\tau)| \ds
  \non
  \\
  && \leq \iotmt \norma{\nabla{\pier\Kt(\mu)(s)}}_2 \norma{\nabla\rho(s)}_4 \norma{\dt\rho(s+\tau)}_4 \ds
  \non
  \\
  && \leq \eps \iot \normaV{\dt\rho(s)}^2 \ds
  + c_\eps {\pier \int_0^{t-\tau} } \normaH{\nabla{\pier\Kt(\mu)(s)}}^2 \normaV{\nabla\rho(s)}^2 \ds
  \non
  \\
  && \leq \eps \intQt |\nabla\dt\rho|^2
  + c \intQt |\dt\rho|^2
  + c_\eps {\pier \int_0^{t-\tau} } \normaH{\nabla{\pier\Kt(\mu)(s)}}^2 \normaV{\nabla\rho(s)}^2 \ds
  \non
  \\
  && \leq \eps \intQt |\nabla\dt\rho|^2 + c 
  + c_\eps {\pier \int_0^{t-\tau} } \normaH{\nabla{\pier\Kt(\mu)(s)}}^2 \normaV{\nabla\rho(s)}^2 \ds.
\Esist
In the last inequality we have used~\eqref{secondastima}.
However, the above estimate has to be improved.
To {\js this end}, we use the regularity theory for linear elliptic equations
and estimates \eqref{pristitau} and \eqref{secondastima} again.
{\pier Indeed, with the help of \eqref{pier1} we have:}
\Bsist
  && \normaV{\nabla\rho(s)}^2
  \leq c \bigl( \normaV{\rho(s)}^2 + \normaH{\Delta\rho(s)}^2 \bigr)
  {\pier {}\leq c\bigl( \normaH{\dt\rho(s)}^2 + 1 \bigr).}
  \non
\Esist
Therefore, the above estimate becomes
\Bsist
  && \intQtmt |\nabla{\pier\Kt(\mu)(s)}| \, |\nabla\rho(s)| \, |\dt\rho(s+\tau)| \ds
  \non
  \\
  && \leq \eps \intQt |\nabla\dt\rho|^2 
  + c_\eps \iot \normaH{\nabla{\pier\Kt(\mu)(s)}}^2 \normaH{\dt\rho(s)}^2 \ds
  + c_\eps \,.
  \label{perquartaBAB}
\Esist
{\pier Analogously, one shows that
\Bsist
  && \intQtmt |\nabla{\pier\Kt(\mu)(s)}| \, |\nabla\rho(s+\tau)| \, |\dt\rho(s+\tau)| \ds
  \non
  \\
  && \leq \eps \intQt |\nabla\dt\rho|^2 
  + c_\eps \iot \normaH{\nabla(\T \Kt(\mu))(s)}^2 \, \normaH{\dt\rho(s)}^2 \ds
  + c_\eps \,. 
  \label{perquartainpiu}
\Esist
}
Thus, by collecting \eqref{perquartaBAA} and \accorpa{perquartaBAB}{perquartainpiu}, we deduce that \eqref{perquartaBA} yields 
\Bsist
  \hskip-2cm && \intQtmt {\pier\kappat(\mu)(s)} \nabla\mu(s) \cdot \nabla \frac 
  {\dt\rho(s+\tau)} {\coeffs} \ds \leq \eps \intQt |\nabla\dt\rho|^2 
  \non
  \\
  \hskip-2cm && \quad {}+  c_\eps \iot {\pier \Big( \normaH{\nabla{\pier\Kt(\mu)(s)}}^2 +   \normaH{\nabla(\T \Kt(\mu))(s)}^2
   \Big)}
\normaH{\dt\rho(s)}^2 \ds
  + c_\eps
  \label{finequartaBA}
\Esist
{\pier for every $\eps>0$.} Let us come to the last term of~\eqref{perquartaB}.
By using the compacness inequality \eqref{compact}, and \eqref{primastima} as well,
we have
\Bsist
  && - \intQtmt \frac {g'(\rho(s)) {\pier g'(\rho(s+\tau)) }}{\coeffs} \, \mu(s) \dt\rho(s) \dt\rho(s+\tau) \ds 
  \non
  \\
  && \leq c \iotmt \norma{\mu(s)}_4 \norma{\dt\rho(s+\tau)}_4 \norma{\dt\rho(s)}_2 \ds
  \non
  \\
  && \leq \eps \iotmt\normaV{\dt\rho(s+\tau)}^2 \ds 
  + c_\eps \iot \norma{\mu(s)}_4^2 \normaH{\dt\rho(s)}^2 \ds
  \non
  \\
  && \leq \eps \iot\normaH{\nabla\dt\rho(s)}^2 \ds + c
  + c_\eps \iot \norma{\mu(s)}_4^2 \normaH{\dt\rho(s)}^2 \ds .
  \label{finequartaBB}
\Esist
Therefore, due to \eqref{finequartaBA} and \eqref{finequartaBB},
\eqref{perquartaB} becomes
\Bsist
  \hskip-1cm&&  \intQt g'(\rho) \dt(\T\mu) \, \dt\rho
  \leq {\pier 2 } \eps \intQt |\nabla\dt\rho|^2 
 + c_\eps \iot \norma{\mu(s)}_4^2 \normaH{\dt\rho(s)}^2 \ds
  \non
  \\
  \hskip-1cm&& \quad
  + c_\eps \iot {\pier \Big( \normaH{\nabla{\pier\Kt(\mu)(s)}}^2 +   \normaH{\nabla(\T \Kt(\mu))(s)}^2
   \Big)} \normaH{\dt\rho(s)}^2 \ds
  + c_\eps \,.
  \label{finequartaB}
\Esist
At this point, we combine \eqref{finequartaZ}, \eqref{finequartaA}, \eqref{finequartaB}
with \eqref{perquarta} and choose $\eps$ small enough.
As the last integral on the \lhs\ is nonnegative since $f_1$ is convex,
we obtain
\Bsist
  && \iO |\dt\rho(t)|^2
  + \intQt |\nabla\dt\rho|^2
  \leq c \iot \phi(s) \normaH{\dt\rho(s)}^2 \ds
  + c
  \non
  \\
  && \quad \hbox{where} \quad
  \phi(s) := {\pier \norma{\mu(s)}_4^2 
+ \normaH{\nabla{\pier\Kt(\mu)(s)}}^2 + \normaH{\nabla(\T \Kt(\mu))(s)}^2
}\, .
  \non
\Esist
As $\phi\in L^1(0,T)$ by \eqref{primastima} and \eqref{daprima},
we can apply the Gronwall lemma and conclude~that
\Beq
  \norma{\dt\rhot}_{\L\infty H\cap\L2V} \leq c.
  \label{quartastima}
\Eeq

\step
Consequence

We have
$-\Delta\rhot+\xit=\psi:=-\dt\rhot+g'(\rhot)\T\mut\in\L\infty H$
due to \eqref{primastima} and \eqref{quartastima}.
Therefore, by a standard argument
(formally multiply by $-\Delta\rhot$ at any fixed time), we deduce that 
both $-\Delta\rhot$ and $\xit$ belong to $\L\infty H$.
Owing to elliptic regularity, we conclude that
\Beq
  \norma\rhot_{\L\infty W} \leq c
  \aand
  \norma\xit_{\L\infty H} \leq c,
  \label{daquarta}
\Eeq
whence also
\Beq
  \norma\rhot_{\LQ\infty}
  + \norma{g(\rhot))}_{\LQ\infty}
  + \norma{g'(\rhot))}_{\LQ\infty}
  + \norma{\pi(\rhot))}_{\LQ\infty}
  \leq c,
  \label{daquartabis}
\Eeq
due to the continuous embedding $W\subset\Linfty$
and the continuity of $g$, $g'$, and~$\pi$ {\pier (note however that $g'$ 
is a bounded function since $g$ is Lipschitz continuous).}

\step 
{\pier Fifth} a priori estimate  

We write \eqref{primatau} as
$\dt\bigl((\coeff)\mu\bigr)=\Delta\Kt(\mu) {\pier {}+{} } g'(\rho)\mu\dt\rho$.
Thus, we have for every $v\in\L2V$
\Bsist
  && \left| \intQ \dt \bigl((\coeff)\mu\bigr) \, v \right|
  = \left| - \intQ \nabla\Kt(\mu) \cdot \nabla v
  {\pier {}+{} } \intQ g'(\rho) \mu \dt\rho \, v \right|
  \non
  \\
  && \leq \norma{\Kt(\mu)}_{{\pier \L2V}} \norma v_{\L2V}
  + \norma{\dt\rho}_{\L\infty H} \norma\mu_{{\pier \L2\Lq}} \norma v_{\L2\Lq}
  \non
  \\
  && \leq \bigl(
    \norma{\Kt(\mu)}_{{\pier \L2V}}
    + c \norma{\dt\rho}_{\L\infty H} \norma\mu_{{\pier \L2\Lq}}
  \bigr) \norma v_{\L2V}.
  \non
\Esist
By accounting for \eqref{primastima}, \eqref{quartastima}, and~{\pier\eqref{daprima}},
we deduce that
\Beq
  \norma{\dt\bigl( (\coefftau)\mut \bigr)}_{\L2\Vp} \leq c.
  \label{sestastima}
\Eeq

\step 
{\pier Sixth} a priori estimate

We test \eqref{primatau} by $\dt\Kt(\mu)=\kappat(\mu)\dt\mu$
and obtain 
\Bsist
  && \intQt (\coeff) \kappat(\mu) |\dt\mu|^2
  + \frac 12 \iO |\nabla\Kt(\mu(t))|^2
  \non
  \\
  && = \frac 12 \iO |\nabla\Kt(\muz)|^2
  - \intQt g'({\pier \rho})\dt\rho\,\mu \dt\Kt(\mu)
  \label{perquinta}
\Esist
{\pier for every $t\in(0,T)$. Note} that the first term on the \lhs\
can be estimated from below as follows
\Beq
  \intQt (\coeff) \kappat(\mu) |\dt\mu|^2
  \geq \intQt \frac {\kappat^2(\mu)} {2\kmax} \, |\dt\mu|^2
  = \frac 1 {2\kmax} \intQt |\dt\Kt(\mu)|^2.
  \label{quintaL}
\Eeq
Now, we deal with the \rhs\ of~\eqref{perquinta}. 
The first term being trivial {\pier thanks to \eqref{hpzero}$_1$}, we come to the second one.
We {\pier invoke} the Young, \holder, and Sobolev inequalities 
and~have 
\Bsist
  && - \intQt g'({\pier \rho})\dt\rho\,\mu \dt\Kt(\mu)
  \leq \frac 1 {4\kmax}  \intQt |\dt\Kt(\mu)|^2
  + c \iot \norma{\mu(s)}_4^2 \norma{\dt\rho(s)}_4^2 \ds 
  \non
  \\
  && \leq \frac 1 {4\kmax}  \intQt |\dt\Kt(\mu)|^2
  + c \iot \norma{\mu(s)}_4^2 \normaV{\dt\rho(s)}^2 \ds .
  \label{quintaR}
\Esist
Now, we observe that \eqref{hpkbis} yields
$\Kt(r)\geq\kmin r-c_*$ for every $r\geq0$, 
where $c_*$ depends on the structural assumptions, only.
Hence, by owing to \eqref{daprimaK} as well, we deduce 
\Bsist
  && \norma{\mu(s)}_4^2 
  \leq c \bigl( \norma{{\pier\Kt(\mu)(s)}}_4^2 + 1 \bigr)
  \leq c \bigl( \normaV{{\pier\Kt(\mu)(s)}}^2 + 1 \bigr)
  \non
  \\
  && \leq c \normaH{\nabla{\pier\Kt(\mu)(s)}}^2 
  + c \normaH{{\pier\Kt(\mu)(s)}}^2 + c 
  \leq c \normaH{\nabla{\pier\Kt(\mu)(s)}}^2 + c 
  \non
\Esist
{\pier for a.a.~$s \in (0,T)$.} 
By combining \eqref{quintaL} and \eqref{quintaR} with \eqref{perquinta},
we obtain
\Bsist
  &&\frac 1 {4\kmax}  \intQt |\dt\Kt(\mu)|^2
  + \frac 12 \iO |\nabla\Kt(\mu(t))|^2
  \leq {\pier {}c +{}} {}c \iot \phi(s) \bigl( \normaH{\nabla{\pier\Kt(\mu)(s)}}^2 + 1 \bigr) \ds
  \non
  \\
  && \quad \hbox{where} \quad
  \phi(s) := \normaV{\dt\rho(s)}^2 .
  \non
\Esist
As $\phi\in L^1(0,T)$ by \eqref{quartastima},
we can apply the Gronwall lemma and conclude that
\Beq
  \norma{\Kt(\mut)}_{\H1H\cap\L\infty V} \leq c.
  \label{quintastima}
\Eeq

\step Consequence

By arguing as we did for \eqref{daprima}, we derive that
\Beq
  \norma\mut_{\L\infty{\Lx6}} \leq c . 
  \label{daquinta}
\Eeq

\step
Limit and conclusion

By the above estimates, there exist a triplet $(\mu,\rho,\xi)$,
with $\mu\geq0$ \aeQ,
and functions $k$ and $\zeta$ such~that
\Bsist
  & \mut \to \mu
  & \quad \hbox{weakly star in $\L\infty{\Lx6}$,}
  \label{convmu}
  \\
  & \rhot \to \rho
  & \quad \hbox{weakly star in $\L\infty W$,}
  \label{convrho}
  \\
  & \dt\rhot \to \dt\rho
  & \quad \hbox{weakly star in $\L\infty H\cap\L2V$,}
  \label{convdtrho}
  \\
  & \xit \to \xi
  & \quad \hbox{weakly star in $\L\infty H$,}
  \label{convxi}
  \\
  & \Kt(\mut) \to k
  & \quad \hbox{weakly star in $\H1H\cap\L\infty V$,}
  \label{convKmu}
  \\
  & \zeta_\tau := (\coefftau)\mut \to \zeta
  & \quad \hbox{weakly star in $\H1\Vp\cap\L\infty{\Lx6}$,}
  \qquad
  \qquad
  \label{convVp}
\Esist
at least for a susequence $\tau=\tau_i{\scriptstyle\searrow}0$.
By \accorpa{convrho}{convdtrho}, \eqref{convKmu},
and the compact embeddings $W\subset\Cx0$ and $V\subset H$,
we can apply well-known strong compactness results
(see, e.g., \cite[Sect.~8, Cor.~4]{Simon})
and have that
\Bsist
  & \rhot \to \rho
  & \quad \hbox{strongly in $\CQ$}
  \label{strongconvrho}
  \\
  & \Kt(\mut) \to k
  & \quad \hbox{strongly in $\C0H$ and \aeQ}.
  \label{strongconvKmu}
\Esist
The weak convergence \eqref{convxi} and \eqref{strongconvrho}
imply that $\xi\in\beta(\rho)$ \aeQ, as is {\pier well known}
(see, e.g., \cite[Prop.~2.5, p.~27]{Brezis}).
The strong convergence \eqref{strongconvrho} also implies 
the Cauchy condition \eqref{cauchyrho}
and that
$\phi(\rhot)\to\phi(\rho)$ strongly in $\CQ$
for every continuous function $\phi:\erre\to\erre$.
We can apply this fact to the functions $g$, $g'$, and~$\pi$
(see \eqref{hpfdueg}).
In~particular, we infer that $\mut g'(\rhot)$ has some weak limit in 
$\L\infty{\Lx6}${\pier : thus, we can identify it with the help of \eqref{convmu}
and} conclude that \eqref{seconda} holds{\pier . Now,} 
we prove that $\mut$ converges to $\mu$ \aeQ.
To this aim, we note that
$\Kt^{-1}$ converges to~$K^{-1}$
uniformly on $[0,R]$ for every $R>0$.
Hence, we see that \eqref{strongconvKmu}
implies $\mut\to K^{-1}(k)$ \aeQ.
By~comparison with~\eqref{convmu}, we deduce that $\mut\to\mu$ \aeQ.
{\pier Next, let us} deal with the subproblem for~$\mu$.
The identity $K^{-1}(k)=\mu$ just proved means that $k=K(\mu)$. 
From the convergence almost everywhere of $\mut$ to $\mu$ we also infer 
that $\zeta_\tau$ converges to $(\coeff)\mu$ \aeQ,
whence $\zeta=(\coeff)\mu$ by comparing with~\eqref{convVp}.
The last term to be identified is the limit of $\eta_\tau:=\mut g'(\rhot)\dt\rhot$.
Precisely, we prove that $\eta_\tau$ converges to $\eta:=\mu g'(\rho)\dt\rho$
weakly in some $L^p$-type space.
We observe that \eqref{convmu} and the convergence almost everywhere of $\mut$
imply that 
\Beq
  \mut \to \mu 
  \quad \hbox{strongly in $\L p{\Lx q}$ for every $p<+\infty$ and $q<6$}
  \label{strongmu}
\Eeq
as is well-known (via the Severini-Egorov theorem).
By~choosing, e.g., $p=q=4$ and combining with 
the weak star convergence of $\dt\rhot$ in $\L\infty H$ (see \eqref{convdtrho})
and the uniform convergence of~$g'(\rhot)$,
we deduce that $\eta_\tau$ converges to $\eta$
weakly in $\L4{\Lx{4/3}}$, thus weakly in $\L2{\Lx{4/3}}$.
At this point, it is \sfw\ to derive~\eqref{wprima}
in an integrated form, namely
\Beq
  \ioT \< \dt \bigl( (\coeff)\mu \bigr)(t) , v(t) > \, dt
  - \intQ \mu g'(\rho)\dt\rho \, v
  + \intQ \nabla K(\mu) \cdot \nabla v
  = 0
  \label{primaint}
\Eeq
for any $v\in\L2V\subset\L2\Lq$,
whence also the {\pier time-}pointwise version \eqref{wprima} itself.
Finally, \eqref{convVp} implies that $\zeta_\tau\to\zeta$
{\pier strongly} in $\C0\Vp$,
thus, $\zeta_\tau(0)\to\zeta(0)$ {\pier strongly} in~$\Vp$,
so that the Cauchy condition \eqref{wcauchymu} is verified as well.
This concludes the proof.


\section{Further properties}
\label{More}
\setcounter{equation}{0}

In this section, we prove Theorems~\ref{Mubdd} and~\ref{Piuregmu}
and make some remarks on the regularity of solutions.
As far as the first result is concerned,
we adapt the argument used in~\cite{CGPS3}.
However, as the first estimate of the proof
has to be derived in a different way,
we prepare a technical lemma.

\Blem
\label{Tecn}
Assume 
$$ {\pier u\in\H1\Vp\cap\L\infty H \quad\hbox{and}\quad u^+\in\H1H\cap\L\infty V ,}$$
and let $\gamma\in\H1V\cap\L\infty W$.
Then, for every $t\in[0,T]$, we have~that
\Beq
  \iot \< \dt u(s) , \gamma(s) u^+(s) > \ds
  = \intQt u \gamma \dt u^+ .
  \label{tecn}
\Eeq
\Elem

\Bdim
As in Remark~\ref{Weakcont}, $u$~is a weakly continuous $H$-valued function
and the pointwise values of $u$ and $u^+$ make sense.
We start from the formula
\Bsist
  && \iot \< \dt u(s) , v(s) > \ds
  = \< u(t) , v(t) >
  - \< u(0) , v(0) >
  - \iot \< u(s) , \dt v(s) > \ds
  \non
  \\
  && = \iO u(t) v(t) 
  - \iO u(0) v(0) 
  - \intQt u \, \dt v,
  \non
\Esist
which holds if $v\in\H1V$.
By an easy \regulariz ation, one sees that it still holds if $v\in\H1H\cap\L2V$.
Now, by applying our assumptions on $u^+$ and $\gamma$
and also owing to the Sobolev inequality~\eqref{sobolev}, 
we have 
\Bsist
  && u^+ \in \L\infty\Lq, \quad
  \nabla u^+ \in \LQ2 , \quad
  \dt u^+ \in \LQ2
  \non
  \\
  && \gamma \in \LQ\infty , \quad
  \nabla\gamma \in \L\infty\Lq, \quad
  \dt\gamma \in \L2\Lq.
\Esist
It follows that $\gamma u^+\in\H1H\cap\L2V$.
Therefore, we obtain
\Bsist
  && \iot \< \dt u(s) , \gamma(s) u^+(s) > \ds
  \non
  \\
  && = \iO u(t) \gamma(t) u^+(t)
  - \iO u(0) \gamma(0) u^+(0)
  - \intQt u \, \dt( \gamma u^+)
  \non
  \\
  && = \iO u(t) \gamma(t) u^+(t)
  - \iO u(0) \gamma(0) u^+(0)
  - \intQt u^+ ( u^+ \dt\gamma + \gamma \dt u^+)
  \non
  \\
  && = \iO \gamma(t) |u^+(t)|^2
  - \iO \gamma(0) |u^+(0)|^2
  - \intQt |u^+|^2 \dt\gamma
  - \intQt \gamma \dt |u^+|^2
  + \intQt u^+ \gamma \dt u^+ 
  \non
  \\
  && = \intQt u \gamma \dt u^+ 
  \non
\Esist
for all $t\in [0,T]$, and the lemma is proved. 
\Edim

\vskip-\bigskipamount
\proofstep
Proof of Theorem~\ref{Mubdd}

Set {\js $\muz^*:=\max\,\{1,\|u_0\|_\infty \}$.}
In the proof performed in~\futuro\cite{CGPS3}
the quantity
$\iO|(\mu(t)-k)^+|^2+\intQt|\nabla(\mu-k)^+|^2$
is estimated for any $k\geq\muz^*$ by testing \eqref{prima}
by $(\mu-k)^+$.
In the present case, the equation to be tested is \eqref{wprima}
instead of~\eqref{prima},
and a~more elaborate procedure is needed.
First of all, we check that $(\mu-k)^+$ is an admissible test function
(this is not obvious since $\nabla\mu$ might not exist in the usual sense).
We recall that $K$ is a strictly increasing mapping from
$[0,+\infty)$ onto itself and that $K^{-1}$
is \Lip\ continuous on the interval $[\smin,+\infty)$,
where $\smin:=K(\rmin)$,
due to~\eqref{hpkbis}.
Therefore, we can choose a strictly increasing map
$\Kz:[0,+\infty)\to[0,+\infty)$ 
that is globally \Lip\ continuous
and coincides with $K^{-1}$ on~$[\smin,+\infty)$.
Hence, we have
$\Kz(K(r))=r$ for every $r\geq\rmin$
and $\Kz(K(r))<\rmin$ for $r<\rmin$.
It follows that $(r-k)^+=(\Kz(K(r))-k)^+$ for every $r\geq0$ if $k\geq\rmin$.
On the other hand, $\Kz(K(\mu))\in\H1H\cap\L2V$ by~\eqref{wregmua}.
Hence, $(\mu-k)^+$ enjoys the same regularity 
and is an admissible test function in~\eqref{wprima}
provided that $k\geq\rmin$.
Thus, we assume $k\geq\max\{\muz^*,\rmin\}$ from now on.
We have from~\eqref{wprima}
\Bsist
  && \iot \< \dt \bigl[ (\coeff)\mu \bigr](s) ,  (\mu(s) - k)^+ > \ds
  + \intQt \nabla K(\mu) \cdot \nabla (\mu-k)^+
  \non
  \\
  && = \intQt \mu \dt g(\rho) \, (\mu-k)^+
  \non
\Esist
for every $t\in[0,T]$,
and a {\pier simple rearrangement} yields
\Bsist
  && \iot \< \dt \bigl[ (\coeff)(\mu-k) \bigr](s) ,  (\mu(s) - k)^+ > \ds
  + \intQt \nabla K(\mu) \cdot \nabla (\mu-k)^+
  \non
  \\
  && = \intQt \dt g(\rho) \, |(\mu-k)^+|^2
  - k \intQt \dt g(\rho)  \, (\mu-k)^+.
  \label{dawprima}
\Esist
Now, noting that $1/(\coeff)\in\H1V\cap\L\infty W$
by \eqref{regrho} and our assumptions on~$g$ (recall \accorpa{hpfg}{hpfdueg}),
we~apply Lemma~\ref{Tecn}
with $u:=(\coeff){\pier (\mu-k)}$ and $\gamma:=1/(\coeff)$
and transform the first term on the \lhs\ as follows:
\Bsist
  && \iot \< \dt \bigl[ (\coeff)(\mu-k) \bigr](s) ,  (\mu(s) - k)^+ > \ds
  = \intQt (\mu-k) \dt \bigl[ (\coeff)(\mu-k)^+ \bigr]
  \non
  \\
  && = \intQt 2 \dt g(\rho) \, |(\mu-k)^+|^2
  + \intQt (\mu-k) (\coeff) \, \dt(\mu-k)^+
  \non
  \\
  && = \frac 12 \intQt \dt \bigl[ (\coeff) |(\mu-k)^+|^2 \bigr]
  + \intQt \dt g(\rho) \, |(\mu-k)^+|^2.
  \non
\Esist
On the other hand, we have 
a.e.\ in the set where $\mu\geq k$
\Beq
  \nabla(\mu-k)^+
  = \nabla\mu
  = \nabla K^{-1}(K(\mu)) 
  = (K^{-1})'(K(\mu)) \nabla K(\mu)
  = \frac 1 {\kamu} \, \nabla K(\mu) .
  \non
\Eeq
Finally, $(\mu(0)-k)^+=0$ since $k\geq\muz^*$.
Therefore, \eqref{dawprima} becomes
\Beq
  \frac 12 \iO (\coefft) |(\mu(t)-k)^+|^2
  + \intQt \kamu |\nabla(\mu-k)^+|^2
  = - k \intQt \dt g(\rho) \, (\mu-k)^+.
  \non
\Eeq
{\pier Since} $g$ is nonnegative and $\kappa(r)\geq\kmin$  for $r\geq k$ 
({\pier because} $k\geq\kmin$), we obtain
\Beq
  \frac 12 \iO |(\mu(t) - k)^+|^2
  + \kmin \intQt |\nabla(\mu-k)^+|^2
  \leq k \intQt |\dt g(\rho)| \, (\mu-k)^+.
  \non
\Eeq
At this point, the argument used in~\futuro\cite{CGPS3}
can be repeated without changes, essentially.
Indeed, the analog of~\eqref{seconda} is never used there,
and the whole proof is based just on the regularity $\dt\rho\in\L\infty H\cap\L2V$. 
In the present case, we have to {\pier exploit} the same regularity 
for~$\dt g(\rho)$
which follows from~\eqref{regrho}.\QED

\Brem
As observed in Remark~\ref{Weakcont},
the component $\mu$ of any weak solution is a weakly continuous $H$-valued function.
We show that 
\Beq
  \mu\in\C0{\Lx p}
  \quad \hbox{for $p\in[1,2)$,\quad {\js and} for $p\in[1,+\infty)$ if $\muz\in\Linfty$}.
  \label{mucont}
\Eeq 
Assume $t_n\to t\in[0,T]$.
We prove that $\mu(t_n)\to\mu(t)$ strongly in $\Lx p$
for $p$ like in~\eqref{mucont}.
{\pier In fact,} $\mu(t_n)$~is bounded in $\Ldue$ in the general case, 
and in $\Linfty$ if $\muz$ is bounded (thanks to Theorem~\ref{Mubdd}).
So, the desired convergence is proved
once we show that $\mu(t_n)\to\mu(t)$ \aeO, at least for a subsequence.
We {\js observe} that \eqref{wregmua} implies $K(\mu)\in\C0H$.
Hence, $K(\mu(t_n))\to K(\mu)$ \aeO, at least for a subsequence.
As $K^{-1}$ is continuous, {\js the claim follows}.
\Erem

\proofstep
Proof of Theorem~\ref{Piuregmu}

As \eqref{rhocont} holds, $g(\rho)$ is continuous and $g'(\rho)$ is bounded.
On the other hand, $\mu$ is bounded too
(by Theorem~\ref{Mubdd} since $\muz$ is bounded). 
Hence, \eqref{prima} can be seen as a linear uniformly parabolic equation for $\mu$
with continuous coefficients and a \rhs\ belonging to~$\L\infty H$.
We have indeed
\Beq
  \dt\mu - \bigl( \coeff \bigr)^{-1} \Delta\mu
  = - \bigl( \coeff \bigr)^{-1} \, \mu \, g'(\rho) \, \dt\rho .
  \non
\Eeq
By owing to $\muz\in W$ and optimal $L^p$-$L^q$-regularity results (see, e.g., \cite[Thm.~2.3]{DHP}), 
we infer that \eqref{piuregmu} holds.\QED

\Brem
\label{Besov}
We notice that the same result {\pier \eqref{piuregmu}} 
holds under an assumption on $\muz$
that is weaker than $\muz\in W$.
The optimal condition involves a proper Besov space
and can give a similar result for a fixed~$p$.
We are going to use \eqref{piuregmu} just with $p=4$
in our proof of the uniqueness of the solution.
It follows that uniqueness still holds for a less regular~$\muz$.
\Erem

\Brem
\label{Ancorareg}
We observe that the case corresponding to an empty interior of $D(\beta)$ is completely trivial.
Indeed, if $D(\beta)=\{r_0\}$, then
$\rho$ takes the constant value~$r_0$, 
$\mu$~solves the corresponding heat equation, and $\xi$ is computed from~\eqref{seconda}.
In the opposite case,
further regularity can be proved under suitable assumptions on the initial data.
For instance, by supposing $\muz$ to be bounded and nonnegative, 
we note that \eqref{seconda} yields
\Beq
  \dt\rho - \Delta\rho + \xi = \mu g'(\rho) - \pi(\rho) \in \LQ\infty .
  \non
\Eeq
So, by assuming that $\inf\rhoz$ and $\sup\rhoz$ belong to the interior of~$D(\beta)$,
one can easily derive that $\xi\in\LQ\infty$.
Indeed, one can formally multiply by $|\xi|^{p-1}\sign\xi$
and estimate $\norma\xi_p$ uniformly with respect to~$p$ 
if this assumption on $\rhoz$ is satisfied.
This implies that 
$\rho\in\W{1,p}{\Lx p}\cap\L p{\Wx{2,p}}$ for every $p<+\infty$
whenever $\rhoz$ is smooth enough.
However, no further regularity can be proved, in general,
since \eqref{seconda} cannot be differentiated,
unless $\beta$ is particular, e.g., like in~\cite{CGPS3}.
By the way, in that case, 
the condition $\xi\in\LQ\infty$ is equivalent to $\inf\rho>0$ and $\sup\rho<1$.
More generally, if $D(\beta)$ is an open interval $(a,b)$ and $\beta$ is a smooth function,
{\js then} $\xi=\beta(\rho)\in\LQ\infty$ implies that $\inf\rho>a$ and $\sup\rho<b$
and a bootstrap technique using both equations can lead to higher regularity.
\Erem


\section{Uniqueness}
\label{Uniqueness}
\setcounter{equation}{0}

In this section, we prove Theorem~\ref{Unicita}.
We observe that the uniqueness of the third component $\xi$ 
follows by comparison in~\eqref{seconda} 
once we prove that the pair $(\mu,\rho)$ is unique.
So, we deal with the first two components, only,
and remind the reader that we can use the further regularity
given by Theorem~\ref{Piuregmu}.
{\js In particular}, by accounting also for~\eqref{regrho} and~\eqref{sobolev}, we~have
\Beq
  |\nabla\mu| \in \L4{\Lx6}
  \aand
  |\nabla\rho| \in \L4{\Lx6} 
  \label{regmurho}
\Eeq
for every solution.
First of all, we rewrite equation~\eqref{prima} in the form
\Beq
  \dt \bigl( \mu / \alpha(\rho) \bigr) - \alpha(\rho) \Delta \mu = 0,
  \label{bravogianni}
\Eeq
where {\js the function} $\alpha:[0,+\infty)\to(0,+\infty)$ is defined by
\Beq
  \alpha(r) := \bigl( \coeffr \bigr)^{-1/2}
  \quad \hbox{for $r\geq 0$} .
  \label{defalpha}
\Eeq
More precisely, we consider the variational formulation of \eqref{bravogianni}
that accounts for the homogeneous Neumann boundary condition
and involves a related unknown function,
namely
\Bsist
  && z := \frac \mu {\alpha(\rho)}
  \aand
  \iO \dt z(t) \, v 
  + \iO \nabla \bigl( \alpha(\rho(t)) \, z(t) \bigr) \cdot \nabla \bigl( \alpha(\rho(t)) v \bigr)
  = 0 
  \non
  \\
  && \quad \hbox{\aat\ and for every $v\in V$}. 
  \label{varbG}
\Esist
{\js Notice that $z$ is bounded since both $\mu$ and $\rho$ are}.
Indeed, \eqref{rhocont} holds and Theorem~\ref{Mubdd} can be applied since $W\subset\Linfty$.
Moreover, $z$~satisfies the analogue of~\eqref{regmurho}.
At this point, we pick two solutions 
$(\mu_i,\rho_i,\xi_i)$, $i=1,2$, and set
\Beq
  a_i := \alpha(\rho_i)
  \aand
  z_i := \mu_i / a_i
  \quad \hbox{for $i=1,2$}
  \non
\Eeq
so that $(z_i,\rho_i)$ satisfy \eqref{varbG}.
In the {\js subsequent} estimates the (varying) value of the constant $\,c\, $
{\js may even depend on the considered solutions},
e.g., through $\norma{z_i}_\infty$.
Our method proceeds as follows.
We write \eqref{varbG} for both solutions and choose
$v=z_1-z_2$ in the difference.
Then we integrate over~$(0,t)$, where $t\in(0,T)$ is arbitrary.
At the same time, we write \eqref{seconda} for both solutions
and multiply the difference by $\rho_1-\rho_2$.
Then we integrate over~$Q_t$.
Finally, we take a suitable linear combination of the {\js resulting} equalities
and perform a number of estimates that lead us to apply the Gronwall lemma.
However, in order to simplify notation
and make the proof more readable, we~set
\Beq
  \mu := \mu_1-\mu_2,\quad
  \rho := \rho_1-\rho_2,\quad
  \xi := \xi_1-\xi_2,\quad
  z := z_1-z_2 ,
  \aand
  a := a_1-a_2, 
  \non
\Eeq
and prepare some auxiliary material before starting.
The next inequalities account for the boundedness and the Lipschitz continuity
of $\alpha$, $\alpha'$, and $1/\alpha$
on the range of~$\rho$ ({\js recall that} $\rho$ is bounded).
We~have
\Bsist
  && |a| = | \alpha(\rho_1) - \alpha(\rho_2)| \leq c |\rho|,
  \non
  \\[1mm]
  && |\nabla a|
  = |\alpha'(\rho_1) \nabla\rho + \bigl( \alpha'(\rho_1) - \alpha'(\rho_2) \bigr) \nabla\rho_2|
  \leq c |\nabla\rho| + c |\nabla\rho_2| \, |\rho|,
  \non
  \\[1mm]
  && |\nabla a_i^{-1}| \leq c |\nabla\rho_i|, 
  \non
  \\[1mm]
  && |\mu|
  \leq |a| \, |z_2| + a_2 |z|
  \leq c |a| + c |z| 
  \leq c |\rho| + c |z| ,
  \non
  \\[1mm] 
  && |\nabla z| = |\nabla \bigl (a_1^{-1} (a_1z) \bigr)|
  \leq c |\nabla(a_1z)| + c |\nabla\rho_1| \, |z|.
  \non
\Esist
{\js In what follows, we will repeatedly use these inequalities} 
without reminding the reader.

\Blem
\label{LemmaA}
We have, for every $t\in[0,T]$,
\Beq
  \intQt |\nabla z|^2 
  \leq c \intQt |\nabla(a_1z)|^2
  + c \iot \bigl( 1 + \norma{\nabla\rho_1(s)}_6^4 \bigr) \, \norma{z(s)}_2^2 \, ds.
  \label{lemmaA}
\Eeq
\Elem

\Bdim
By the preliminary inequalities just {\js stated}, we have
\Beq
  \intQt |\nabla z|^2
  \leq c \intQt |\nabla(a_1z)|^2
  + c \intQt |\nabla\rho_1|^2 \, |z|^2.
  \label{perlemmaA}
\Eeq
Now, by using the \holder\ and Sobolev inequalities (see~\eqref{sobolev}), 
we~obtain
\Bsist
  && \intQt |\nabla\rho_1|^2 \, |z|^2
  \leq \iot \norma{\nabla\rho_1(s)}_6^2 \, \norma{z(s)}_6 \, \norma{z(s)}_2 \, ds
  \non
  \\
  && \leq c \iot \norma{\nabla\rho_1(s)}_6^2 \,\bigl( \norma{\nabla z(s)}_2 + \norma{z(s)}_2 \bigr)\, \norma{z(s)}_2 \, ds
  \non
  \\
  && \leq c \iot \norma{\nabla\rho_1(s)}_6^2 \, \norma{z(s)}_2^2 \, ds
  \non
  \\
  && \quad {}
  + \eps \intQt |\nabla z|^2
  + c_\eps \iot \norma{\nabla\rho_1(s)}_6^4 \, \norma{z(s)}_2^2 \, ds
  \non
  \\
  && \leq \eps \intQt |\nabla z|^2
  + c_\eps \iot \bigl( 1 + \norma{\nabla\rho_1(s)}_6^4 \bigr) \, \norma{z(s)}_2^2 \, ds,
  \non
\Esist
where $\eps>0$ is arbitrary.
Hence, \eqref{lemmaA} follows by combining this with~\eqref{perlemmaA}
and then choosing $\eps$ small enough.
\Edim

\Blem
\label{LemmaB}
Let $k\in\L4{\Lx6}$. Then we have 
\Bsist
  && \intQt k^2 (|z|^2 + |\rho|^2)
  \leq \eps \intQt
    \bigl( |\nabla(a_1z)|^2 + |\nabla\rho|^2
    \bigr)
  \non
  \\
  && \quad {}
  + c_\eps \iot
    \bigl( 1 + \norma{\nabla\rho_1(s)}_6^4 + \norma{k(s)}_6^4 \bigr)
    \bigl( \norma{z(s)}_2^2 + \norma{\rho(s)}_2^2 \bigr) \, ds 
  \label{lemmaB}
\Esist
for every $\eps>0$ and every $t\in[0,T]$.
\Elem

\Bdim
By the \holder\ and Sobolev inequalities (see~\eqref{sobolev}), we have
\Bsist
  && \intQt k^2 (|z|^2 + |\rho|^2)
  \leq \iot \norma{k(s)}_6^2 \,
    \bigl(
      \norma{z(s)}_6 \, \norma{z(s)}_2
      + \norma{\rho(s)}_6 \, \norma{\rho(s)}_2
    \bigr) \, ds
  \non
  \\
  && \leq c  \iot \norma{k(s)}_6^2
    \bigl(
      \norma{\nabla z(s)}_2 \, \norma{z(s)}_2
      + \norma{z(s)}_2^2
      + \norma{\nabla\rho(s)}_2 \, \norma{\rho(s)}_2
      + \norma{\rho(s)}_2^2
    \bigr) \, ds
  \non
  \\
  && \leq \eps \iot 
    \bigl(
      \norma{\nabla z(s)}_2^2 
      + \norma{\nabla\rho(s)}_2^2 
    \bigr) \, ds
  + c_\eps \iot 
    \bigl(
      \norma{k(s)}_6^4
      + \norma{k(s)}_6^2
    \bigr)
    \bigl(
      \norma{z(s)}_2^2 
      + \norma{\rho(s)}_2^2 
    \bigr) \, ds.
    \non
\Esist
By applying Lemma \ref{LemmaA} and denoting 
the constant that appears in~\eqref{lemmaA} by~$C$,
we can continue and obtain
\Bsist
  && \intQt k^2 (|z|^2 + |\rho|^2)
  \non
  \\
  && \leq \eps
    \tonde{
      C \intQt |\nabla(a_1z)|^2 
      + C \iot \bigl( 1 + \norma{\nabla\rho_1(s)}_6^4 \bigr) \, \norma{z(s)}_2^2 \, ds
      + \intQt |\nabla\rho|^2
    }
  \non
  \\
  && \quad {}
  + c_\eps \iot 
    \bigl(
      1 + \norma{k(s)}_6^4
    \bigr)
    \bigl(
      \norma{z(s)}_2^2 + \norma{\rho(s)}_2^2 
    \bigr) \, ds.
    \non
  \non
\Esist
Hence, \eqref{lemmaB} immediately follows.
\Edim

At this point, we can start with our program.
However, in order to make the argument more transparent,
we deal with the first equation only, for a while.
We have
\Beq
  \frac 12 \iO |z(t)|^2
  + \intQt \bigl(
    \nabla(a_1z_1) \cdot \nabla(a_1z)
    - \nabla(a_2z_2) \cdot \nabla(a_2z)
  \bigr) = 0 .
  \non
\Eeq
It is convenient to transform the last integrand as follows:
\Bsist
  && \nabla(a_1z_1) \cdot \nabla(a_1z)
  - \nabla(a_2z_2) \cdot \nabla(a_2z)
  \non
  \\
  && = |\nabla(a_1z)|^2
  + \nabla(a_1z_2) \cdot \nabla(a_1z)
  - \nabla(a_2z_2) \cdot \nabla(a_1z)
  + \nabla(a_2z_2) \cdot \nabla(az)
  \non
  \\
  && = |\nabla(a_1z)|^2
  + \nabla(az_2) \cdot \nabla(a_1z)
  + \nabla\mu_2 \cdot \nabla(az).
  \non
\Esist
Then, the above equality becomes
\Beq
  \frac 12 \iO |z(t)|^2
  + \intQt |\nabla(a_1z)|^2
  = - \intQt \nabla(az_2) \cdot \nabla(a_1z)
  - \intQt  \nabla\mu_2 \cdot \nabla(az),
  \label{diffprima} 
\Eeq
and we estimate each term of the {\rhs} separately.
We immediately have
\Bsist
  && - \intQt \nabla(az_2) \cdot \nabla(a_1z)
  \leq \frac 14 \intQt |\nabla(a_1z)|^2
  + 2 \intQt \bigr( z_2^2 |\nabla a|^2 + a^2 |\nabla z_2|^2 \bigr)
  \non
  \\
  && \leq \frac 14 \intQt |\nabla(a_1z)|^2
  + c \intQt \bigl( |\nabla\rho|^2 + |\nabla\rho_2|^2 \, |\rho|^2 \bigr)
  + c \intQt |\nabla z_2|^2 \, |\rho|^2  
  \non
  \\
  && \leq \frac 14 \intQt |\nabla(a_1z)|^2
  + C_1 \intQt |\nabla\rho|^2 
  + c \intQt \bigl( |\nabla\rho_2|^2 + |\nabla z_2|^2 \bigr)\, |\rho|^2,  
  \label{rhsA}
\Esist
where we have marked the constant {\js that} we want to refer to by terming it~$C_1$.
We treat the last term of \eqref{diffprima} as follows:
\Bsist
  && - \intQt \nabla\mu_2 \cdot \nabla(az)
  \leq \intQt |\nabla\mu_2| \bigl( |a| \, |\nabla z| + |z| \, |\nabla a| \bigr)
  \non
  \\
  && \leq c \intQt |\nabla\mu_2| \bigl( 
    |\nabla(a_1z)| \, |\rho|
    + |z| \, |\nabla\rho_1| \, |\rho|
    + |z| \, |\nabla\rho|
    + |z| \, |\nabla\rho_2| \, |\rho|
  \bigr)
  \non
  \\
  && \leq \frac 14 \intQt |\nabla(a_1z)|^2
  + c \intQt |\nabla\mu_2|^2 \, |\rho|^2
  + \intQt |\nabla\rho|^2 
  \non
  \\
  && \quad {}
  + c \intQt |\nabla\mu_2|^2 \, |z|^2
  + c \intQt \bigl( |\nabla\rho_1|^2 + |\nabla\rho_2|^2 \bigr) \, |\rho|^2 .
  \label{rhsB}
\Esist
Now, we deal with the second equation.
Testing the difference of \eqref{seconda} by $\rho$
as mentioned at the beginning, easily yields
\Beq
  \frac 12 \iO |\rho(t)|^2 
  + \intQt |\nabla\rho|^2
  + \intQt \xi \rho
  = \intQt \bigl( \mu_1 g'(\rho_1) - \mu_2 g'(\rho_2) - \pi(\rho_1) + \pi(\rho_2) \bigr) \rho\,.
  \label{diffseconda}
\Eeq
{\js We note} that the last integral on the \lhs\ of \eqref{diffseconda}
is nonnegative by monotonicity,
while the integrand on the \rhs\ can be estimated as follows:
\Bsist
   &&{\pier \bigl( \mu_1 g'(\rho_1) - \mu_2 g'(\rho_2) - 
   \pi(\rho_1) + \pi(\rho_2) \bigr) \rho }
  \non
  \\[1mm]
  && \leq \bigl(
    |\mu| \, |g'(\rho_1)|
    + |\mu_2| \, |g'(\rho_1) - g'(\rho_2)|
    + |\pi(\rho_1) - \pi(\rho_2)|
  \bigr) |\rho|
  \non
  \\[1mm]
  && \leq |g'(\rho_1)| \, |\mu| \, |\rho|
  + c |\mu_2| \, |\rho|^2
  + c |\rho|^2
  \leq c \bigl( |\mu|^2 + |\rho|^2 \bigr) 
  \leq c \bigl( |z|^2 + |\rho|^2 \bigr) .
  \non
\Esist
Thus, by {\js inspecting} the coefficients of the integral $\intQt|\nabla\rho|^2$
that appear on the \rhs s of \eqref{rhsA} and~\eqref{rhsB},
it is clear that it is convenient to multiply \eqref{diffseconda}
by $\,C_1+2\,$ before {\js adding it} to~\eqref{diffprima}.
Once such a care is taken,
it is \sfw\ to deduce that
\Bsist
  && \iO |z(t)|^2
  + \intQt |\nabla(a_1z)|^2
  + \iO |\rho(t)|^2
  + \intQt |\nabla\rho|^2
  \non
  \\
  && \leq c \intQt \bigl( |\nabla\mu_2|^2 + 1 \bigr)\, |z|^2  
  + c \intQt \bigl( |\nabla\mu_2|^2 + |\nabla\rho_1|^2 + |\nabla\rho_2|^2 + |\nabla z_2|^2 + 1 \bigr) \, |\rho|^2 
  \non
  \\
  && \leq c \intQt \bigl( |\nabla\mu_2| + |\nabla\rho_1| + |\nabla\rho_2| + |\nabla z_2| + 1 \bigr)^2 \,
   \bigl( |z|^2 + |\rho|^2 \bigr) . 
  \non
\Esist
At this point, 
we recall that \eqref{regmurho} hold for $\mu_i$, $\rho_i$, and~$z_i$,
so that we can apply Lemma~\ref{LemmaB} with {\js $k=|\nabla\mu_2|+|\nabla\rho_1|+ |\nabla\rho_2| + |\nabla z_2|+1$}.
Then, we choose {\js $\eps>0$} small enough and use the Gronwall lemma.
We conclude that $z=0$ and $\rho=0$,
whence $(\mu_1,\rho_1)=(\mu_2,\rho_2)$ {\js follows}.

{\pier
\section*{Acknowledgements}

PC and JS gratefully acknowledge the warm hospitality of the WIAS 
in Berlin and the IMATI of CNR in Pavia. Some financial support 
came from the MIUR-PRIN Grant 2008ZKHAHN {\sl ``Phase transitions, 
hysteresis and multiscaling''} for PC and GG, and from the 
FP7-IDEAS-ERC-StG Grant \#200947 (BioSMA) for PC and JS.
The work of JS was also supported by the DFG Research Center 
{\sc Matheon} in Berlin.
} 


\vspace{3truemm}

\Begin{thebibliography}{10}

\bibitem{Barbu}
V. Barbu,
``Nonlinear semigroups and differential equations in Banach spaces'',
Noord\-hoff,
Leyden,
1976.

\bibitem{bcdgss}
E. Bonetti, P. Colli, W. Dreyer, G. Gilardi, G. Schimperna,
J. Sprekels, On a model for phase separation in binary
alloys driven by mechanical effects, {\it Phys.~D} {\bf 165} (2002), 48-65. 

\bibitem{bds}
{E. Bonetti, W. Dreyer, G. Schimperna, Global solutions to 
a generalized {C}ahn-{H}illiard equation with viscosity,
{\it Adv. Differential Equations} {\bf 8} (2003) 231-256.}

\bibitem{Brezis}
H. Brezis,
``Op\'erateurs maximaux monotones et semi-groupes de contractions
dans les espaces de Hilbert'',
North-Holland Math. Stud.
{\bf 5}, North-Holland,
Amsterdam, 1973.

{\bibitem{CN} B.D. Coleman, W. Noll, 
The thermodynamics of elastic materials with heat conduction and viscosity, 
{\it Arch. Rational Mech. Anal.} {\bf 13} (1963) 167-178.} 

\bibitem{CGPS1} 
P. Colli, G. Gilardi, P. Podio-Guidugli, J. Sprekels,
Existence and uniqueness of a global-in-time solution
to a phase segregation problem of the Allen-Cahn type, 
{\it Math. Models Methods Appl. Sci.} {\bf 20} (2010)
519-541.

\bibitem{CGPS2} 
P. Colli, G. Gilardi, P. Podio-Guidugli, J. Sprekels,
A temperature-dependent phase segregation problem 
of the Allen-Cahn type, {\it Adv. Math. Sci. Appl.} 
{\bf 20} (2010) {219-234.}

\bibitem{CGPS3} 
P. Colli, G. Gilardi, P. Podio-Guidugli, J. Sprekels,
Well-posedness and long-time behaviour 
for a nonstandard viscous Cahn-Hilliard 
system, {\it SIAM J. Appl. Math.} {\bf 71} (2011) 1849-1870.

\bibitem{CGPS4} 
P. Colli, G. Gilardi, P. Podio-Guidugli, J. Sprekels,
An asymptotic analysis for a nonstandard Cahn-Hilliard system with viscosity, 
preprint WIAS-Berlin n.~1652 (2011), pp.~1-17. {\js Accepted for publication
in {\it  Discrete Contin. Dyn. Syst. Ser. S.}}  

\bibitem{CGPS5} 
P. Colli, G. Gilardi, P. Podio-Guidugli, J. Sprekels,
Distributed optimal control of a nonstandard 
system of phase field equations, 
{\pier {\it Contin. Mech. Thermodyn.}  
doi:10.1007/s00161-011-0215-8}

\bibitem{CGPS7} 
P. Colli, G. Gilardi, P. Podio-Guidugli, J. Sprekels,
Global existence for a strongly coupled 
Cahn-Hilliard system with viscosity, 
preprint WIAS-Berlin n.~1691 (2012), pp.~1-15. 
{\pier Accepted for publication in 
{\it Boll. Unione Mat. Ital. (9).}} 

\bibitem{CGS1} 
P. Colli, G. Gilardi, J. Sprekels,
Analysis and optimal boundary 
control of a nonstandard system of phase field equations, 
{\it Milan J. Math.} doi:10.1007/s00032-012-0181-z

\bibitem{DHP}
R. Denk, M. Hieber, J. Pr\"uss,
Optimal $L^p$-$L^q$-estimates for parabolic 
boundary value problems with inhomogeneous data.
{\it Math. Z.} {\bf 257} (2007) 193-224.

\bibitem{Fremond}
M. Fr\'emond,
``Non-smooth Thermomechanics'',
Springer-Verlag, Berlin, 2002.

{\bibitem{FG} 
E. Fried, M.E. Gurtin, 
Continuum theory of thermally induced phase transitions based on an order 
parameter, {\it Phys. D} {\bf 68} (1993) 326-343.}

\bibitem{GPS} 
M. Grasselli, H. Petzeltov\'a, and G. Schimperna, Asymptotic behavior of a nonisothermal viscous Cahn-Hilliard equation with inertial term. {\it J. Differential Equations} {\bf 239} (2007) 38-60. 

\bibitem{Gurtin} 
M. Gurtin, Generalized Ginzburg-Landau and
Cahn-Hilliard equations based on a microforce balance,
{\it Phys.~D\/} {\bf 92} (1996) 178-192.

\bibitem{mr} 
{A. Miranville, A Rougirel, Local and asymptotic analysis 
of the flow generated by the
{C}ahn-{H}illiard-{G}urtin equations,  
{\it Z. Angew. Math. Phys.} {\bf 57} (2006) 244-268.}

\bibitem{Podio}
P. Podio-Guidugli, 
Models of phase segregation and diffusion of atomic species on a lattice,
{\it Ric. Mat.} {\bf 55} (2006) 105-118.

\bibitem{Ros} 
{R. Rossi, On two classes of generalized viscous {C}ahn-{H}illiard
equations, {\it Commun. Pure Appl. Anal.} {\bf 4} (2005) 405-430.}

\bibitem{Simon}
J. Simon,
{Compact sets in the space $L^p(0,T; B)$},
{\it Ann. Mat. Pura Appl.} {\bf 146} (1987) 65-96.

\End{thebibliography}

\End{document}

\bye